\documentclass[12pt]{article}
{\begin{Sbox}\begin{minipage}}%
{\end{minipage}\end{Sbox}\fbox{\TheSbox}}%
\usepackage{booktabs}
\usepackage{graphicx}
\usepackage{float}
\usepackage{subfigure}
\usepackage{sidecap}
\usepackage[psamsfonts]{amssymb}
\usepackage{amsmath, amsthm, anysize, enumerate, color, url}
\usepackage[authoryear]{natbib}
 \usepackage{epstopdf}%%\makeatletter \renewcommand\@biblabel[8]{${#8}$}\makeatother
\usepackage{longtable}
\usepackage{bm}
\usepackage[colorlinks,citecolor=blue,urlcolor=blue]{hyperref}
\usepackage{amsmath}
\usepackage{mathrsfs}
\usepackage{geometry}
\usepackage{caption}
\usepackage{subfigure}
\usepackage{float}
\usepackage{breakurl}

\numberwithin{equation}{section}%
\makeatletter
\newcommand\Hsub{\@startsection{subsection}{2}%
 {0pt}{-\baselineskip}{.2\baselineskip}%
 {\normalsize\upshape}}
\makeatother
\newtheorem{theorem}{Theorem}[section] %如果不采用章节号做前缀，则不用[section]
\newtheorem{definition}{Definition}[section]
\newtheorem{lemma}{Lemma}[section]
\newtheorem{remark}{Remark}[section]
\newtheorem{proposition}{Proposition}[section]

\newcommand{\pozhehao}{\kern0.3ex\rule[0.8ex]{1.5em}{0.095ex}\kern0.3ex}

\usepackage{setspace}

\begin{document}

\vspace*{0.15cm}
\begin{center}
{\Large\bf Asymptotic Theory for Differentially Private Generalized $\beta$-models
with Parameters Increasing}
\vskip 1.2\baselineskip
{\large Yifan Fan$^{1}$, \ \ Huiming Zhang$^{^2*}$, \ \ Ting Yan$^{1}$}
\vskip 0.2\baselineskip
{\sl Department of Statistics, Central China Normal University, Wuhan, 430079, China$^{1}$}
\vskip 0.1\baselineskip
\vskip 0.2\baselineskip
{\sl School of Mathematical Sciences and Center for Statistical Science, Peking University, Beijing, 100871, China$^{2}$}
\vskip 0.1\baselineskip

\end{center}
\footnotetext[2]{$^*$Correspondence author.\\
Email addresses: fanyifan${\_}7@$163.com (Y. Fan), zhanghuiming$@$pku.edu.cn (H. Zhang ),\\ tingyanty$@$mail.ccnu.edu.cn (T. Yan)}.
\vskip 1.5mm

\begin{abstract}
Modelling edge weights play a crucial role in the analysis of network data, which reveals the extent of relationships among individuals. Due to the diversity of weight information, sharing these data has become a complicated challenge in a privacy-preserving way. In this paper, we consider the case of the non-denoising process to achieve the trade-off between privacy and weight information in the generalized $\beta$-model.
Under the edge differential privacy with a discrete Laplace mechanism, the Z-estimators from estimating equations for the model parameters are shown to be consistent and asymptotically normally distributed. The simulations and a real data example are given to further support the theoretical results.


{\noindent } 	
\end{abstract}

\vskip 5 pt

\section{Introduction}
\par
With the rapid development of computer and network technology,  the analysis of network data  has aroused widespread concerns in various fields.
Unfortunately, collecting, storing, analyzing and sharing these data is challenging, due to the privacy of individuals (e.g.,  financial transactions).
Besides, more privacy protection may reduce the validity of data [\cite{Duncan: Keller-McNulty:Stokes:2004}].
Many approaches have been proposed to guarantee the trade-off between individual privacy and the utility of published data, which focus on data encryption, identity authentication, data perturbation  [\cite{Samarati and Sweeney1998, Fung:Wang:Yu:2007, Machanavajjhala: Gehrke:Kifer D:Venkitasubramaniam:2006, Ghinita:Tao:Kalnis:2008, Li:Li:Venkatasubramanian:2007,Aggarwal and Yu 2007}].   \cite{Dwork2006a} proposed a rigorous notion of privacy named $\varepsilon$-differential privacy to control strong worst-case privacy risks. More formally, adding or removing a single record in the dataset does not have a serious effect on the outcome of any analysis. Starting from \cite{Dwork2006a}, various types of data and queries were widely applied by researchers under differential privacy constraints [\cite{Holohan:Leith:Mason:2017, McSherry :Talwar:2007,Wasserman:Zhou:2010}].

 Random graphs are powerful statistical tools in the study of network data. These graph models are based on degree sequences $d$, which are used in modelling the realistic networks.
 In the undirected case, the $\beta$-model is well-known for the binary network, renamed by \cite{Chatterjee:Diaconis:Sly:2011}. Many scholars have focused on the studies of the $\beta$-model [\cite{Jackson.(2008), Lauritzen.(2008), Blitzstein and Diaconis(2011)}]. \cite{Chatterjee:Diaconis:Sly:2011} proved the existence and consistency of the maximum likelihood estimator (MLE) of the $\beta$-model as the number of parameters goes to infinity. \cite{Yan:Xu:2013} further derived its asymptotic normality. On the other hand,  edge weights reveal the strength of relationships among individuals, which are critical for understanding many phenomena. For example, in friendship networks, we can assign close friends with a higher weight and acquaintances or normal friends with a lower weight, which are also referred to as the strong tie and weak tie reported by \cite{Granovetter(1993)}.
To this point, \cite{Hillar:Wibisono:2013} studied the maximum entropy distributions on weighted graphs with the $\beta$--model as the special case and
proved the consistency of the MLE under the assumption that all parameters are bounded by a constant;
 \cite{Yan:Zhao:Qin:2015} proved the asymptotic normality of the MLE.

 In the privacy analysis of network data, the raw data is published via pre-processing so that the confidential and sensitive information is captured as less as possible.
One of the popular approaches is to add some noises $e$  into the degree sequence $d$.
 For example, \cite{Hay:Li:Miklau:Jensen:2009} applied the Laplace noise-addition
mechanism  to release the degree partition of a graph, and designed to reduce the error with the $\ell_2$--norm between the true and released degree partitions.
However, the process of adding noises is often ignored when summary statistics are published in a privacy-preserving way.
As a result, the estimated parameters  may not be consistent, even not exist [e.g., \cite{Hay:Li:Miklau:Jensen:2009}].
\cite{Duchi2018} illustrated that the estimator operated on private data has a larger error than the non-private estimator.
Based on privatized data, estimating summary statistics and estimating parameters of models are totally different problems.
To this point, \cite{Karwa2016} paid attention to the noise addition process through the denoised method  to achieve valid inference and
obtained the consistency and asymptotic normality of a differential privacy estimator in the $\beta$-model.

In this paper, we adopt the non-denoised method to establish the asymptotic properties of the Z-estimator of the parameter in the generalized $\beta$-model with finite discrete weighted edges under the discrete Laplace mechanism,
which is different from the work of \cite{Karwa2016}. Moreover, \cite{Karwa2016} only considered binary edges.
In some scenarios, edge weights play important roles in the analysis of network data.
For example, weighted social networks often provide a more realistic representation of the complex social interactions among individuals than binary networks [e.g., \cite{Farine(2014)}].
Furthermore,  edge weights may further increase the risk of privacy disclosure, due to  the diversity of weight-related information.
For instance, edge weights represent the numbers of co-written papers in a coauthorship network. A hacker can easily identify  an author  via the total number of published papers [\cite{Li:Shen:Lang:Dong:2016}].
In the generalized $\beta$-model, each node is assigned one parameter, so the number of parameters increases with $n$.
The asymptotic properties for the increasing dimensional $Z$-estimator cannot directly be followed from the classical $Z$-estimation theory; see chapter 5 of \cite{Vaart1998}.
Therefore, based on \cite{Yan:Xu:2013}, we alternatively show that the $Z$-estimator of the parameter involving the noisy degree sequence is asymptotically consistent and normally distributed in the generalized $\beta$-model under edge differential privacy constraints.

\par The organization of this paper is as follows. In Section \ref{s2}, we first introduce some notations and definitions of our results. Subsequently, we obtain the asymptotic normality of the $Z$-estimator in the generalized $\beta$-model involving noisy degree sequence $\bar{\textbf{d}}=\textbf{d}+\textbf{e}$, where $\textbf{d}$ is the sufficient statistic and $ \textbf{e}$  are some noises from the discrete Laplace distribution. In Section \ref{s3}, we give some simulation results to support our theories. We further present a data example application, which is from a community of $27$ Grevy's zebras. A summary follows in Section \ref{s4}. All proofs are contained in Section \ref{s5}.

\section{Main Results}\label{s2}

\subsection{Notations }
\par For a vector $\textbf{x}=(x_{1},\cdots, x_{n})^{T}\in R^{n},$  the $\ell_{\infty}$-norm of $\textbf{x}$ is denoted by $\Vert \textbf{x}\Vert_{\infty}=\mathop {\max }\limits_{1 \le i \le n} |{x_i}|$. For an $n\times n$ matrix $J=(J_{ij})$, $\Vert J \Vert_{\infty}$ denotes the matrix norm induced by the $\Vert\cdot \Vert_{\infty}$-norm on vectors in $R^{n}:$
$$\Vert J \Vert_{\infty}=\max_{\textbf{x}\neq 0}\cfrac{\Vert J\textbf{x}\Vert_{\infty}}{\Vert \textbf{x}\Vert_{\infty}}=\max_{1\leq i\leq n}\sum^{n}_{j=1}\vert J_{ij} \vert,$$
i.e., the maximum absolute row sum norm.

We define another matrix norm
$\|\cdot\|$ for a matrix $A=\left(a_{i, j}\right)$ by $\|A\| :=\max _{i, j}\left|a_{i, j}\right|$, and let $\Vert \textbf{x}\Vert_{1}=\sum_{i}\vert x_{i}\vert$ be the  $\ell_1$-norm for a general vector $\textbf{x}.$ We say that $a_{n}=\Omega(r_n)$ if there exists a real constant $c>0$ and there exists an integer constant $n_{0}\geq 1$ such that $a_{n}\geq c r_n$ for every $n\geq n_{0}$.
\par Let $D$ be an open convex subset of $R^{n}.$ An $n\times n$ function matrix $G(\textbf{x})$ whose elements $G_{ij}(\textbf{x})$ are functions on a vector $\textbf{x},$ is Lipschitz continuous w.r.t the max norm on $D$ if there exists a real number $\kappa$ such that for any $\textbf{v}\in R^{n}$ and any $\textbf{x},\textbf{y}\in D,$
$$\Vert G(\textbf{x})(\textbf{v})-G(\textbf{y})(\textbf{v})\Vert_{\infty}\leq \kappa\Vert\textbf{x}-\textbf{y} \Vert_{\infty}\Vert\textbf{v} \Vert_{\infty}, $$
where $\kappa$ may depend on $n$ but it is independent of $\textbf{x}$ and $\textbf{y}$. For every fixed $n$, $\kappa$ is a constant.

\subsection{Edge Differential Privacy}
\par In the contexts of network data, there are two main variants of differential privacy: edge differential privacy (EDP)[\cite{Nissim:Raskhodnikova:Smith:2007}] and node differential privacy (NDP) [\cite{Hay:Li:Miklau:Jensen:2009}, \cite{Kasiviswanathan2013}], which are based on the different definitions of graph neighbors. Specifically, EDP guarantees that
released databases do not reveal the addition or removal of a special edge,  while NDP hides the addition or removal of a node (along with its edges) in a graph $G.$ In this paper, we refer to EDP, where two graphs $G$ and $G'$ are said to be neighbors if  they differ in exactly one edge.
\begin{definition}
{\rm (Edge differential privacy). Let $\varepsilon\geq 0$ be a privacy parameter. A randomized mechanism (or a family of conditional probability distributions) $\mathcal{Q}(\cdot|G)$ is $\varepsilon$- edge differentially private if
$$ \sup\limits_{G,G'\in \mathcal{G},\delta(G,G')=1}\sup_{S\in\mathcal{S}}\cfrac{\mathcal{Q}(S|G)}{\mathcal{Q}(S|G')}\leq e^{\varepsilon},$$
where $\mathcal{G}$ is the set of all undirected graphs of interest on $n$ nodes, $\delta(G,G')$ is the number of edges on which $G$ and $G'$ differ, $\mathcal{S}$ is the set of all possible outputs (or the support of  $\mathcal{Q}$ ).
}\end{definition}
\par The above definition of EDP is based on ratios of probabilities. Generally, the data curator chooses an appropriate privacy parameter  $\varepsilon$ to achieve the trade-off between privacy and validity.
 As the value of $\varepsilon$ is extremely small, more privacy is protected. Under EDP constraints, changing one record in the dataset cannot affect seriously on the distribution of the output. For example, a hospital can release some medical information about their patients to the public, while simultaneously ensuring very high levels of privacy in the case of EDP. This is because EDP offers a guarantee no matter whether or not the patient participates in the study, the probability of a possible output is almost the same. As a result, an attacker can not find whether a single individual  is in the original database or not.  As we know, the effective implementation of $\varepsilon$-differential privacy is associated with the magnitude of additional random noise. To this end, \cite{Dwork:McSherry:Nissim:Smith:2006b} introduced the notion of global sensitivity, which is referred to as the maximum $\ell_1$- norm among various dataset pairs $(G,G').$
\begin{definition}
{\rm(Global sensitivity). Let $f:\mathcal{G}\rightarrow \mathcal{R}^{k}.$ The global sensitivity of $f$ is defined as
$$\triangle_G(f)=\max_{\delta(G,G')=1}\Vert f(G)-f(G')\Vert_{1},$$
where $\Vert \cdot\Vert_{1}$ is the $\ell_1$-norm for vector.
}\end{definition}
\par Although there are many mechanisms for releasing the output of any function $f$ under differential privacy, the Laplace mechanism is the most common one. \cite{Karwa2016} presented a discrete Laplace mechanism to achieve edge differential privacy, which is given below.

 Let $f:\mathcal{G}\rightarrow \mathcal{Z}^{k},$ and let $Z_{1},\ldots, Z_{k}$ be independent and identically distributed discrete Laplace random variables with probability mass function defined by
$$P(Z=z)=\frac{1-\lambda}{1+\lambda}\lambda^{|z|}, z\in \mathcal{Z}, \lambda\in(0,1).$$
Then the algorithm which outputs $f(G)+(Z_{1},\ldots, Z_{k})$ with inputs $G$ is $\varepsilon$-edge differentially private, where $\varepsilon=-\bigtriangleup_G(f)\log\lambda.$

\par Based on the definition of differential privacy, \cite{Dwork:McSherry:Nissim:Smith:2006b} found that any function of a differentially private mechanism is also differentially private, as follow: Let $f$ be an output of an $\varepsilon$-differentially private mechanism and $g$ be any function. Then $g(f(G))$ is also $\varepsilon$-differentially private. This result indicates that any post-processing done on the noisy degree sequence obtained as an output of a differentially private mechanism is also differentially private.
 \par More generally, we may consider the \emph{skew discrete Laplace mechanism}. When the positive noises and negative noises arising with different probability law, the skew discrete Laplace distribution [\cite{Kozubowski2006}] as a discretization of non-symmetric Laplace distribution could be used.
The skew Laplace distribution is useful in applications to communications, engineering, and finance and economics, see \cite{Kotz2012} and references therein. For more information on the \emph{skew discrete Laplace mechanism}, see the supplementary material for details.

\subsection{Estimation}
\par   Let $\mathcal{G}_{n}$ be a simple undirected graph including $n$ nodes. Let $a_{ij}$ be the weight of edge $(i,j)$ , $1\leq i \neq j \leq n$, taking values from the set $\{0,1,\ldots,q-1\}$. Let $A=(a_{ij})$ be the adjacency matrix of  $\mathcal{G}_{n}$. Note that $\mathcal{G}_{n}$ has no self-loops, $a_{ii}=0.$ Define $d_{i}=\sum_{j\neq i}a_{ij}$ and $\bm{d}=(d_{1},\cdots,d_{n})^{T}$ as the degree sequence of  $\mathcal{G}_{n}$. The density or probability mass function on $\mathcal{G}_{n}$ with respect to some canonical measure $\nu$ has the exponential-family random graph models with the degree sequence as sufficient statistic, i.e.,
\begin{align*}
p(\mathcal{G}_{n};\bm{\alpha})=\exp(\bm{\alpha}^{T}\bm{d}-Z(\bm{\alpha})),
\end{align*}
where $Z(\bm{\alpha})$ is the normalizing constant, $\bm{\alpha}=(\alpha_{1},\cdots,\alpha_{n})^{T}$ is a vector parameter.

We assume that the edge weights $\{a_{ij}\}$ are independently multinomial random variables with the  probability mass function:
\begin{align}\label{eq:ERG}
P(a_{ij}=a)=\frac{e^{a(\alpha_{i}+\alpha_{j})}}{\sum_{k=0}^{q-1}e^{k(\alpha_{i}+\alpha_{j})}}, \ \ \ \ a=0,1,\ldots, q-1.
\end{align}
where $q\ge 2$ is a fixed number of the class. Thus the likelihood  of $\bm{\alpha}$ is
\[L(\alpha ) \propto \prod\limits_{j\neq i} {\prod\limits_{a = 0}^{q - 1} {{{[P{\rm{(}}{a_{ij}} = a{\rm{)]}}}^{1({a_{ij}} = a)}}} } ,\]
which gives the log-likelihood of $\bm{\alpha}$,
\begin{align*}
\log L(\bm{\alpha}) &\propto \sum\limits_{j\neq i} {\sum\limits_{a = 0}^{q - 1} {\left[ {1({a_{ij}} = a) \cdot \left\{ {a({\alpha _i} + {\alpha _j}) - \log \left( {\sum\limits_{k = 0}^{q - 1} {{e^{k({\alpha _i} + {\alpha _j})}}} } \right)} \right\}} \right]} }.
\end{align*}
So the log-partition function
 in \eqref{eq:ERG} is $Z(\bm{\alpha})=\sum_{j\neq i}\log \sum_{k=0}^{q-1}e^{(\alpha_{i}+\alpha_{j})k}.$ This model is  a direct generalization of the $\beta$-model, which only considers the dichotomous edges.

Moreover, the first order condition for the log-likelihood function w.r.t. ${\alpha _i}$ are
\[\frac{{\partial \log L(\bm{\alpha})}}{{\partial {\alpha _i}}} =\sum_{j\neq i}^{n} {{a_{ij}} - \sum\limits_{j=1;j\neq i} {\sum\limits_{a = 0}^{q - 1} {\frac{{a{e^{a({\alpha _i} + {\alpha _j})}}}}{{\sum\limits_{k = 0}^{q - 1} {{e^{k({\alpha _i} + {\alpha _j})}}} }}} } } ,i = 1,2, \ldots ,n.\]
Let  $e_{i}$ be a noise independently drawn from discrete Laplace distribution with parameter $\lambda_n$. We output $\bar{d_{i}}:=d_{i}+e_{i}$ using the discrete Laplace mechanism. However, the degree $ d_{i}=\sum_{j\neq i}^{n}a_{ij}$ of vertex $i$ is not attainable, since the observed degree contains unknown noise in private date set. We resort to the moment equations which are given by the following system of functions:
\[
    \begin{split}
&F_{i}(\bm{\alpha})=\bar{d_{i}}-E(\bar{d_{i}})=\bar{d_{i}}-E(d_{i}), \ \ \ \ i=1,\ldots,n,\\
&F(\bm{\alpha})=(F_{1}(\bm{\alpha}), \ldots,F_{n}(\bm{\alpha}))^{T}.
\end{split}
    \]

Under this case, since adding or removing an edge can change the degree of at most two nodes, by $1$ each, the global sensitivity for the degree sequence $d$ is $2$. Therefore, we have the privacy parameter
$$\varepsilon_n:=-\bigtriangleup_G(f)\log\lambda_n=-2 \log\lambda_n.$$
So, ${\lambda _n} = \exp ( - \frac{{{\varepsilon _n}}}{2})$.

 We use $\hat{\bm{\alpha}}$ to denote the Z-estimator of $\bm{\alpha}$ satisfying $F(\hat{\bm{\alpha}})=0.$
 Since the noises $e_{i}$'s ($i=1, 2, \cdots, n$) are independently drawn from symmetric discrete Laplace distribution with parameter $\lambda_{n}$,  $E(e_{i})=0$.
Note that $d_i$ is a sum of edge weights $a_{ij}$'s ($j=1, \ldots n, j\neq i$).
So we have
\[
    \begin{split}
    E(\bar{d}_{i})&=E(d_{i}+e_{i})=E(d_{i})=\sum^{n}\limits_{j=1;j\neq i}E(a_{ij})\\
&=\sum^{n}\limits_{j=1;j\neq i}\sum_{k=0}^{q-1} kP(a_{ij}=k)=\sum^{n}\limits_{j=1;j\neq i}\sum^{q-1}_{k=0}\frac{ke^{k(\alpha_{i}+\alpha_{j})}}{\sum^{q-1}\limits_{k=0}e^{k(\alpha_{i}+\alpha_{j})}}.
\end{split}
    \]
Therefore the moment-based estimating equations with noisy degree sequences are
\begin{equation}\label{eq1}
\bar{d_{i}}=\sum_{j=1;j\neq i}^{n}\sum_{a=0}^{q-1}\frac{ae^{a(\hat{\alpha}_{i}+\hat{\alpha}_{j})}}{\sum_{k=0}^{q-1}e^{k(\hat{\alpha}_{i}+\hat{\alpha}_{j})}}, \ \ \ \ i=1,\ldots, n.
    \end{equation}

\subsection{Consistency and Asymptotical Normality }
\par In this section, we obtain that the Z-estimator of the parameter involving noisy degree sequence is asymptotically consistent and normally distributed.

\par Given $m,M>0,$ we say an $n\times n$ matrix $V_{n}=(v_{ij})$ belongs to the matrix class $\mathcal{L}_{n}(m,M)$ if $V_{n}$ is a symmetric nonnegative matrix satisfying
$$v_{ii}=\sum^{n}_{j=1,j\neq i}v_{ij};\\\\ M\geq v_{ij}=v_{ji}\geq m>0, i\neq j.$$
Generally, the inverse of $V_n$, $V^{-1}_{n}$, does not have a closed form. \cite{Yan:Xu:2013} proposed a simple matrix $\bar{S}=(\bar{s}_{ij})$ to approximate $V^{-1}_{n},$ where $\bar{s}_{ij}=\cfrac{\delta_{ij}}{v_{ii}}-\cfrac{1}{v_{..}},$ $\delta_{ij}$ is the Kronecker delta function, and $v_{..}=\sum_{i,j}^{n}(1-\delta_{ij})v_{ij}=\sum^{n}_{i=1}v_{ii}.$

\par Similar to \cite{Yan:Qin:Wang:2016b}, let the parameter vector $\boldsymbol{\alpha}=(\alpha_{1},\cdots,\alpha_{n})^{T}$ belong to the  symmetric parameter space
$$D=\{\bm{\alpha}\in R^{n}: -Q_{n}\leq \alpha_{i}+ \alpha_{j}\leq Q_{n}, 1\leq i< j\leq n\}, $$
 where $\{ Q_{n}\}$ is the sequence of upper bound of the parameters. Let $F'(\bm{\alpha})$ be the
Jacobian matrix of $F(\bm{\alpha})$ at $\bm{\alpha}$, then for $i,j=1,\ldots n,$
\[
    \begin{split}
&\frac{\partial F_{i}}{\partial\alpha_{i}}=\sum_{j=1;j\neq i}^{n}\frac{\sum_{0\leq k<l\leq q-1}(k-l)^2e^{(k+l)(\alpha_{i}+\alpha_{j})}}{(\sum_{a=0}^{q-1}e^{a(\alpha_{i}+\alpha_{j})})^{2}},\\
&\frac{\partial F_{i}}{\partial\alpha_{j}}=\frac{\sum_{0\leq k<l\leq q-1}(k-l)^2e^{(k+l)(\alpha_{i}+\alpha_{j})}}{(\sum_{a=0}^{q-1}e^{a(\alpha_{i}+\alpha_{j})})^{2}},\ \ j=1,\ldots n; j\neq i.
\end{split}
    \]
Here we need assume $V_n:=F'(\bm\alpha)\in \mathcal{L}_{n}(m,M),$ i.e., $v_{ii}=\cfrac{\partial F_{i}}{\partial\alpha_{i}}$ and $v_{ij}=\cfrac{\partial F_{i}}{\partial\alpha_{j}}, $ where
 \begin{equation}
  \label{eq2}m=(2(1+e^{^{Q_{n}}}))^{-1} \text{~and~}\ M=\cfrac{q^{2}}{2}.
 \end{equation}

First, we give the convergence rate of the $\ell_{\infty}$ error which directly shows the consistency of the parameters under some mild conditions.

\begin{theorem}\label{th21}
{\rm Consider the discrete Laplace mechanism with ${\lambda _n} = \exp ( -\epsilon_n/2)$ and assume that $\bm{\alpha}\in D$ and $e^{Q_{n}}=o((n/\log n)^{\frac{1}{12}})$, where $D=\{\bm{\alpha}\in R^{n}: -Q_{n}\leq \alpha_{i}+\alpha_{j}\leq Q_{n}, \text{for}\ 1\leq i<j\leq n \}$. If $\varepsilon_{n}\geq c\sqrt{\log n}$ (denoted as $\varepsilon_{n}=\Omega(\sqrt{\log n})$), where $c\geq 4$ is a constant,   then as $n$ goes to infinity, $\bm{\hat{\alpha}}$ exists and satisfies
\begin{equation}\label{eq3}
\Vert\bm{\hat{\alpha}}-\bm{\alpha} \Vert_{\infty}=O_{p}(e^{3Q_{n}}\sqrt{\frac{\log n}{n}})=o_{p}(1).
    \end{equation}}
\end{theorem}
\begin{remark}
{\rm  In Theorem \ref{th21},  we use  Newton's method to obtain the existence and consistency of $\bm{\hat{\alpha}}$. This indicates that the Z-estimator of the parameter $\bm{\alpha}$ involving a noisy sequence is accurate under the non-denoised process. If ${Q_{n}}$ is bounded and thus $\bm{\alpha}$ is a sparse vector, this convergence rate matches the oracle inequality $\| {\bm{\hat \beta}  - \bm{\beta ^*}} \|_{\infty} = {O_p}(\sqrt {\frac{{\log p}}{n}} )$ for the Lasso estimator in the linear model with $p=n-1$ -dimensional true parameter vector $\bm{\beta ^*}$ and the sample size $n$, see \cite{Lounici08}.}
\end{remark}
Second, we get the asymptotic normality of the estimator in the restricted parameter space under the slower rate condition for $e^{Q_{n}}$ compared with the rate for $e^{Q_{n}}$ in Theorem \ref{th21}, as follow.

\begin{theorem}\label{th22}
{\rm Consider the discrete Laplace mechanism with ${\lambda _n} = \exp ( -\epsilon_n/2)$ and assume that the conditions in Theorem \ref{th21} hold. If we assign a smaller $e^{Q_{n}}=o\left(n^{1/18}/(\log n)^{1/9}\right)$, then as $n$ goes to infinity, for any fixed $r\geq1,$ the vector
$$(v_{11}^{1/2}(\hat{\alpha}_{1}-\alpha_{1}),\cdots,v_{rr}^{1/2}(\hat{\alpha}_{r}-\alpha_{r}) )\stackrel {d} { \rightarrow } N_r( \mathbf { 0 } , \mathbf {I}_r )$$
where $\mathbf {I}_r$  is a $r \times r$ identity matrix.}
\end{theorem}

\begin{remark}
{\rm By Theorem \ref{th22}, for any fixed $i$, the convergence rate of $\hat{\alpha}_{i}$ is $1/ v^{1/2}_{ii},$ when $\varepsilon_{n}=\Omega(\sqrt{\log n})$. Since $(n-1)(2(1+e^{^{Q_{n}}}))^{-1} \leq v_{ii}\leq (n-1)q^2/2,$  the convergence rate is between $O(n^{-1/2}e^{Q_n/2})$ and $O(n^{-1/2})$ , which is  the same as the non-privacy estimator (\cite{Yan:Qin:Wang:2016b}). }
\end{remark}

The proofs of Theorems \ref{th21} and \ref{th22} are postponed in the Appendix section. After deriving the theoretical results, numerical studies are carried out in the next section to verify the asymptotic properties of the $Z$-estimate. Theorem \ref{th22} can be also used to construct a confidence interval for the parameters. For instance, an approximate $1-\alpha$ confidence interval for $\alpha_{i}-\alpha_{j}$ is $\hat{\alpha}_{i}-\hat{\alpha}_{j}\pm Z_{1-\alpha/2}(1/\hat{v}_{ii}+1/\hat{v}_{jj})^{1/2},$ where $Z_{1-\alpha/2}$ is the $1-\alpha$-quantile of the standard normal distribution, $\hat{v}_{ii}$ and $\hat{v}_{jj}$ are the $Z$-estimates of $v_{ii}$ and $v_{jj}$ by replacing all $\alpha_{i}$  with their $Z$-estimates.

\section{Numerical studies}\label{s3}
\subsection{Simulations}

We first consider simulations under a discrete weight $q=3$.
In this case, we evaluate asymptotic properties by simulating finite sample data in finite networks. We consider the changes of $n, \varepsilon$ and $L$. Based on \cite{Yan2015} and \cite{Yan:Leng:Zhu:2016a}, the setting of the true parameter vector $\bm{\alpha^{*}}$ takes a linear form. Specifically, we set $\alpha^{*}_{i}=(n-i+1)L/n,$ for $i=1,\cdots,n$. We discuss three distinct values for $L$, $L=0, \log(\log n), (\log n)^{1/2}$, respectively. We simulate three distinct values for $\varepsilon:$ one is fixed $(\varepsilon=2)$ and the other two values tend to zero with $n,$ i.e., $\varepsilon=\log(n)/n^{1/4},\log(n)/n^{1/2}.$ Here we discuss three values for $n$, $n=100$, $200$ and $500$. Each simulation is repeated $10,000$ times.

\begin{figure}[!htb]
\centering
\subfigure[finite discrete weights ($q=3$)]{
\setlength{\abovecaptionskip}{0.cm}
\setlength{\belowcaptionskip}{-0.cm}
\includegraphics[height=6in, width=6in]{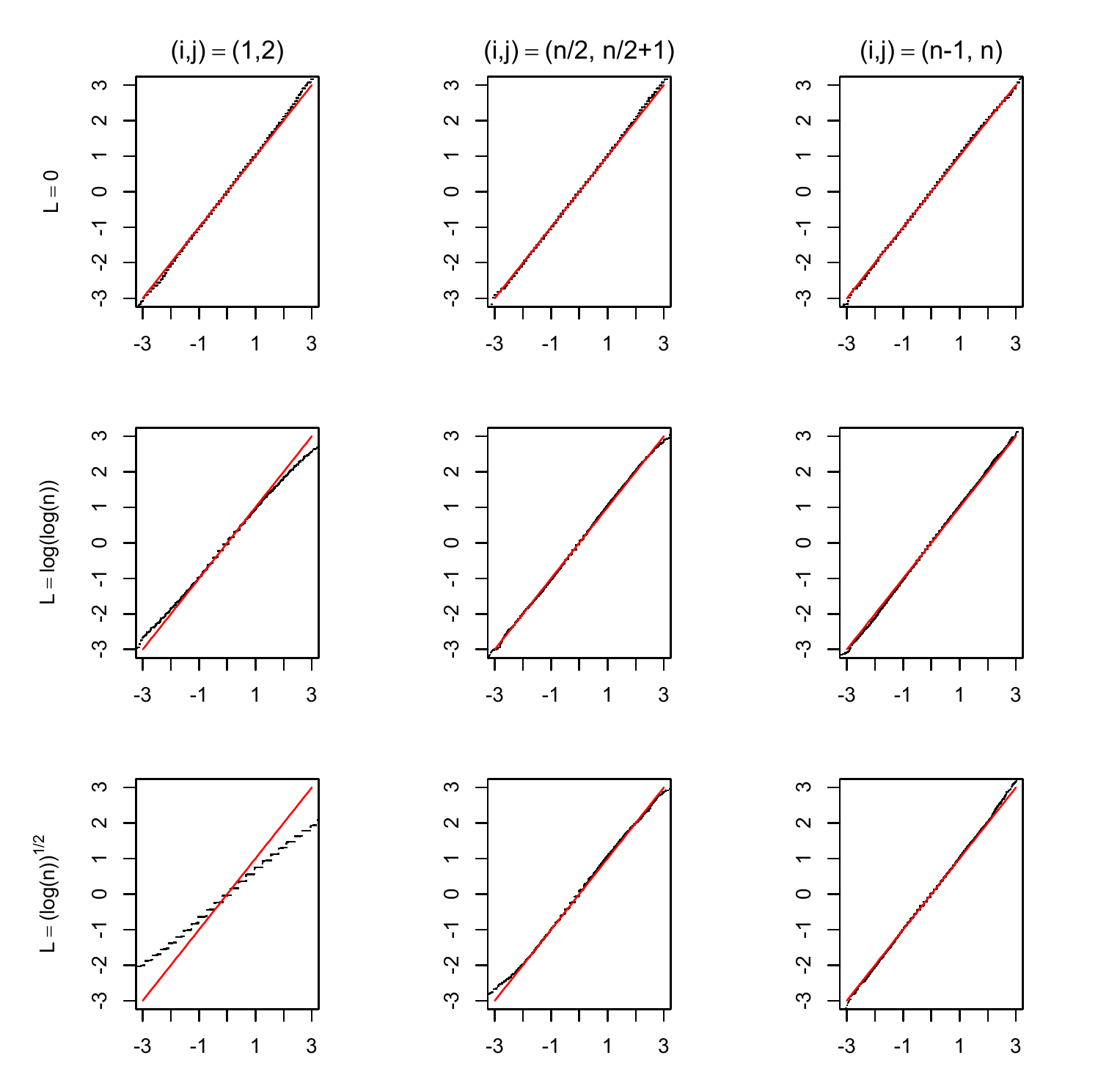}
}
\caption{The QQ plots ($n=100,\epsilon=2$) .}
\label{graph1}
\end{figure}

\par By Theorem \ref{th22}, $\hat{\xi}_{ij}=[\hat{\alpha}_{i}-\hat{\alpha}_{j}-(\alpha_{i}^{*}-\alpha_{j}^{*})]/(1/\hat{v}_{ii}+1/\hat{v}_{jj})^{1/2}$ converges to the standard normal distribution, where $\hat{v}_{ii}$ is the estimator of $v_{ii}$ by replacing $\alpha_{i}^{*}$ with $\hat{\alpha}_{i}.$ Hence, we apply the quantile-quantile (QQ) plot to demonstrate the asymptotic normality of $\hat{\xi}_{ij}.$ Three special pairs $(1,2), (n/2,n/2+1)$ and $(n-1, n)$ for $(i,j)$ are presented in Figure \ref{graph1}. Further, we list the coverage probability of the $95\%$ confidence interval, the length of the confidence interval, and the frequency that the estimate does not exist.
\par For $\varepsilon=2,\log(n)/n^{1/4}$, the QQ-plots under $n=100, 200$ and $500$ are similar. Thus, we here only show the QQ-plots for $\hat{\xi}_{ij}$ under the case of $\varepsilon=2$ and $n=100$ in Figure \ref{graph1} to save space. In Figure \ref{graph1}, the horizontal and vertical axes are the theoretical and empirical quantiles, respectively, and the red lines correspond to the reference lines $y=x.$ From Figure \ref{graph1}, we see that for fixed pair $(i,j)=(1,2)$, the empirical quantiles coincide well with the ones of the standard normality for noisy estimates (i.e., $\hat{\xi}_{ij}$) expect for $L=(\log n)^{1/2}$. When $L=(\log n)^{1/2}$, notable deviations exist for pair $(1,2)$ in Figure \ref{graph1}. For other pairs $(n/2,n/2+1)$ and $(n-1,n)$, the approximation of asymptotic normality is good when $L=0, \log(\log n), (\log n)^{1/2}.$  While

\begin{table}
\centering
\caption{Estimated coverage probabilities of $\alpha_{i}^{*}-\alpha_{j}^{*}$  for pair $(i,j)$ as well as the length of confidence intervals (in square brackets), and the probabilities that the estimate does not exist (in parentheses).}
\scalebox{0.8}{
\begin{tabular*}{13cm}{lllll}
\hline
$n$\ \ \ & $(i,j)$ \ \ \ & $L=0$\ \ \  &$L=\log(\log n)$ \ \ \ &$L=(\log n)^{1/2}$\\
\hline
&&&$\epsilon=2$\\
\hline
$100$ &$(1,2)$&$94.63[0.35](0)$ &$96.75[0.81](0.41) $&$99.75[1.13](31.79)$\\
&$(50,51)$&$94.80[0.35](0)$&$94.79[0.55](0.41)$&$95.07[0.73](31.79)$\\
&$(99,100)$&$94.90[0.35](0)$&$94.04[0.41](0.41)$&$94.55[0.46](31.79)$\\\\
$200$&$(1,2)$&$94.40[0.25](0)$ &$97.68[0.62](0)$ &$99.86[0.85](6.42)$\\
&$(50,51)$&$94.51[0.25](0)$&$94.54[0.41](0)$&$95.74[0.54](6.42)$\\
&$(99,100)$&$94.44[0.16](0)$&$94.45[0.30](0)$&$94.92[0.33](6.42)$\\\\
$500 $&$(1,2)$&$95.17[0.16](0)$ &$98.62[0.42](0)$ &$99.98[0.57](0.03)$\\
&$(50,51)$&$94.79[0.16](0)$&$94.67[0.27](0)$&$96.81[0.36](0.03)$\\
&$(99,100)$&$95.05[0.16](0)$&$94.91[0.19](0)$&$95.20[0.21](0.03)$\\
\hline
&&&$\epsilon=\log (n)/n^{1/4}$\\
\hline
$100$ &$(1,2)$&$94.41[0.35](0)$ &$95.80[0.81](1.42)$ &$99.48[1.12](51.83)$\\
&$(50,51)$&$94.40[0.35](0)$&$94.07[0.56](1.42)$&$93.96[0.73](51.83)$\\
&$(99,100)$&$94.39[0.35](0)$&$93.66[0.41](1.42)$&$94.12[0.46](51.83)$\\\\
$200 $&$(1,2)$&$94.40[0.25](0)$ &$96.90[0.62](0.03)$ &$99.70[0.84](17.92)$\\
&$(50,51)$&$94.34[0.25](0)$&$94.10[0.41](0.03)$&$95.13[0.54](17.92)$\\
&$(99,100)$&$94.24[0.25](0)$&$94.09[0.30](0.03)$&$94.42[0.33](17.92)$\\\\
$500 $&(1,2)&$95.03[0.16](0)$ &$98.30[0.42](0)$ &$99.94[0.57](0.73)$\\
&$(50,51)$&$94.60[0.16](0)$&$94.53[0.27](0)$&$96.40[0.36](0.73)$\\
&$(99,100)$&$94.92[0.16](0)$&$94.82[0.19](0)$&$95.00[0.21](0.73)$\\
\hline
&&&$\epsilon=\log (n)/n^{1/2}$\\
\hline
$100$ &$(1,2)$&$88.65[0.35](0) $&$83.92[0.82](61.94)$ &$87.50[1.03](99.60)$\\
&$(50,51)$&$88.84[0.35](0)$&$81.56[0.57](61.94)$&$77.50[0.76](99.60)$\\
&$(99,100)$&$88.00[0.35](0)$&$84.73[0.42]((61.94)$&$67.50[0.47](99.60)$\\\\
$200$ &$(1,2)$&$90.03[0.25](0)$ &$81.25[0.63](40.95)$ &$86.21[0.78](99.71)$\\
&$(50,51)$&$89.57[0.25](0)$&$81.69[0.42](40.95)$&$72.41[0.53](99.71)$\\
&$(99,100)$&$89.27[0.25](0)$&$87.33[0.30](40.95)$&$93.10[0.33](99.71)$\\\\
$500 $&$(1,2)$&$91.52[0.16](0)$ &$83.64[0.43](12.78)$ &$93.55[0.56](99.38)$\\
&$(50,51)$&$91.38[0.16](0)$&$84.41[0.28](12.78)$&$75.81[0.37](99.38)$\\
&$(99,100)$&$91.59[0.16](0)$&$89.46[0.19](12.78)$&$87.10[0.21](99.38)$\\
\hline
\end{tabular*}}
\label{tab1}
\end{table}
\noindent  $\varepsilon=\log(n)/n^{1/2}$, the approximation of asymptotic normality respecting to $\hat{\xi}_{ij}$ is bad, see Figure 1 of the supplementary material.
\par The coverage probability of the $95\%$ confidence interval for $\alpha^{*}_{i}-\alpha^{*}_{j}$, the length of the confidence interval, and the frequency that the estimate does not exist, are reported in Table \ref{tab1}. The length of the confidence interval is related to  $L$ and $n$. That is, the length increases as $L$ increases, or the length decreases as $n$ increases. Under the case of $\varepsilon=2,\log(n)/n^{1/4}$, the coverage frequencies of pair $(1,2)$ are higher than the nominal level $95\%$ expect for $L=0$; for other pairs $(n/2,n/2+1)$ and $(n-1,n)$, the coverage frequencies are all close to the nominal level $95\%$ for all $L$, where the ones are the closest at $n = 500$. For $\varepsilon=\log(n)/n^{1/2}$, the coverage frequencies are  lower than the nominal level $95\%$ for all $L$. This indicates that as $\varepsilon$ reduces to a specific value (e.g., $\log(n)/n^{1/2}$), notable deviations exist between the coverage frequencies and the nominal level $95\%$, especially  the probabilities of the non-existent estimates are very high when $L=(\log n)^{1/2}$.

Second, we compare the simulation results between with the denoising process [\cite{Karwa2016}] and without the denoising process  in the case of $q=2$.
 Here the settings of $\hat{\bm{\alpha}}^{*}, \bar{\bm{\alpha}}^{*}, L$ and $\varepsilon$ are the same
as those in the first simulation. Here, we only consider $n=100, 200$ without $500$. Each simulation is repeated $10,000$ times.

\begin{figure}[!htb]
\centering
\subfigure[finite discrete weights ($q=2$)]{
\setlength{\abovecaptionskip}{0.cm}
\setlength{\belowcaptionskip}{-0.cm}
\includegraphics[height=6in, width=6in]{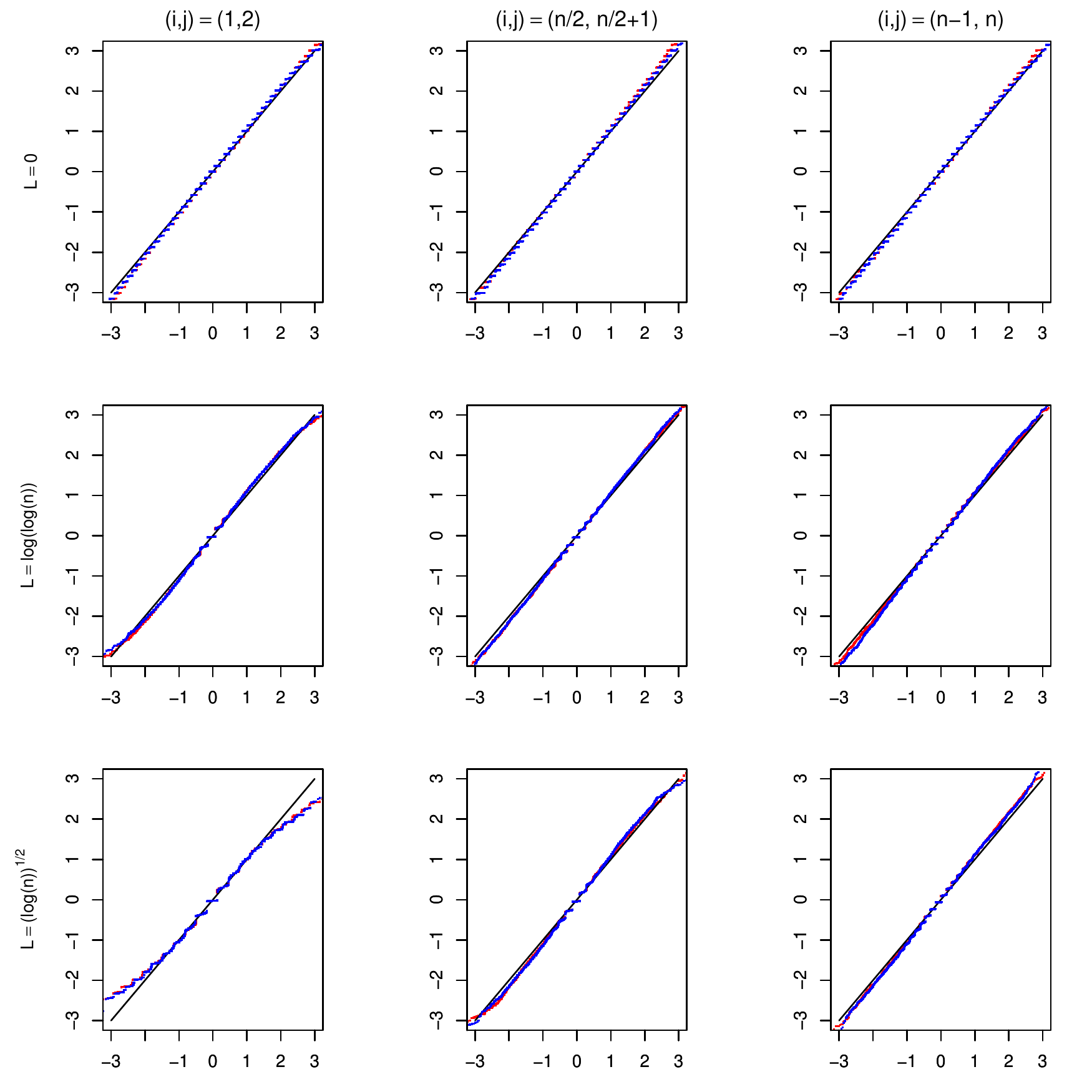}
}
\caption{The QQ plots of the non-denoised and denoised estimates ($n=100,\epsilon=2$) .}
\label{graph2}
\end{figure}

According to the results in  \cite{Karwa2016}, $\bar{\xi}_{ij}=[\bar{\alpha}_{i}-\bar{\alpha}_{j}-(\alpha_{i}^{*}-\alpha_{j}^{*})]/(1/\bar{v}_{ii}+1/\bar{v}_{jj})^{1/2}$ converges to the standard normal distributions,
where $\bar{\alpha}_{i}$ is the estimate of $\alpha_{i}$ with the denoising process and  $\bar{v}_{ii}$ is  the estimate of $v_{ii}$ by replacing ${\alpha}_{i}$ with $\bar{\alpha}_{i}$.
We apply the quantile-quantile (QQ) plot  and record
the coverage probability of the 95$\%$ confidence interval, the length of confidence interval, and the frequency that the estimate does not exist,
to compare the performance of $\hat{\xi}_{ij}$ and $\bar{\xi}_{ij}$.
The QQ-plots are shown in Figure \ref{graph2} and numeric comparison results are given in Table \ref{table2}.
In Figure \ref{graph2}, the QQ-plots for both $\hat{\xi}_{ij}$ denoted by the red  color and $\bar{\xi}_{ij}$ denoted by the blue color are very close
and  coincide well with the ones of the standard normality when $\varepsilon=2,\log(n)/n^{1/4}$ and $L\le \log(\log n)$. (We only show the QQ-plots of $\varepsilon=2$ and $n=100$ in Figure \ref{graph2} to save space and the other cases are similar.)
This indicates that the parameter estimates are nearly the same with and without the denoising process.
However, when $\varepsilon=\log(n)/n^{1/2}$, the approximation of asymptotic normality of both $\hat{\xi}_{ij}$ and $\bar{\xi}_{ij}$ is not good, see Figure 2 of the supplementary material.

In Table \ref{table2},  Type ``A" and ``B" represent the estimates without and with the denoised process, respectively. From this table, we can see that the difference between both estimates is very small.
Similar to the analysis of Table \ref{tab1}, the length of confidence interval increases as $L$ increases and decreases as $n$ increases. Under the case of $\varepsilon=2,\log(n)/n^{1/4}$, the coverage frequencies of all pairs are all close to the nominal level $95\%$ when $L=0, \log(\log n);$ for $L=(\log n)^{1/2}$, both the non-denoised and denoised estimates often failed to exist for $n=100$, while $n = 200$ the non-existent frequencies of estimates are lower. For $\varepsilon=\log(n)/n^{1/2}$, the coverage frequencies for both non-denoised and denoised estimates exist a great gap compared with the nominal level $95\%$ for all $L$, and the probabilities of the non-existent  estimates also increase as $L$ increases.

\begin{table}[h!]\centering
\footnotesize
\label{tab2}
\caption{Estimated coverage probabilities of $\alpha_{i}^{*}-\alpha_{j}^{*}$  for pair $(i,j)$ as well as the length of confidence intervals (in square brackets), and the probabilities that the estimate does not exist (in parentheses). Type ``A" denotes the estimate with the non-denoised process and  ``B" denotes the estimate with the denoised process. }
\begin{tabular}{lcccccc}
 \hline
  $n$ & $(i,j)$ & Type & $L=0$ & $L=\log(\log(n))$ & $L=(\log n)^{1/2}$   \\
  \hline
  \multicolumn{7}{c}{$\epsilon =2 $}\\
   \hline
 $100$ & (1,2)  & A &93.62[0.57](0)&93.38[1.01](1.25)&$97.38[1.46](43.88)$\\

       &        & B &93.79[0.57](0)&93.76[1.01](1.30)&97.35[1.47](44.17)\\\\

       &(50,51) &A &93.44[0.57](0)&93.57[0.76](1.25)&93.16[0.94](43.88)\\

       &        & B &93.62[0.57](0)&93.48[0.76](1.30)&93.19[0.94](44.17)\\\\

       &(99,100)&A &93.42[0.57](0)&93.63[0.63](1.25)&93.16[0.68](43.88)\\

       &       & B &93.82[0.57](0)&92.90[0.63](1.30)&93.43[0.68](44.17)\\\\
  \hline
   $200$ & (1,2)&A &94.55[0.40](0)&93.80[0.75](0.03)&96.60[1.11](9.94)\\
       &        &B &94.85[0.40](0)&94.04[0.75](0.03)&96.75[1.11](10.55)\\\\
       &(100,101)&A &95.09[0.40](0)&94.28[0.55](0.03)&93.93[0.68](9.94)\\
       &        &B &94.77[0.40](0)&94.58[0.55](0.03)&93.87[0.68](10.55)\\\\
       &(199,200)&A &94.89[0.40](0)&94.20[0.45](0.03)&93.75[0.48](9.94)\\
       &       &B&94.28[0.40](0)&94.23[0.45](0.03)&93.56[0.48](10.55)\\
  \hline
  \multicolumn{7}{c}{$\epsilon =\log(n)/n^{1/4} $}\\
   \hline
 $100$ & (1,2)  & A &92.62[0.58](0)&91.31[1.02](4.46)&96.04[1.46](65.58)\\
       &        & B &92.74[0.58](0)&91.74[1.02](5.14)&95.99[1.45](66.34)\\\\
       &(50,51) &A &92.56[0.58](0)&91.88[0.76](4.46)&91.11[0.95](65.58)\\
       &        &B &92.67[0.58](0)&92.00[0.76](5.14)&91.21[0.95](66.34)\\\\
       &(99,100)&A &92.70[0.58](0)&92.58[0.63](4.46)&91.81[0.68](65.58)\\
       &       & B &92.78[0.58](0)&91.79[0.64](5.14)&92.13[0.68](66.34)\\\\
\hline
$200$ & (1,2)  & A &94.14[0.40](0)&92.03[0.76](0.19)&95.34[1.12](26.08)\\
       &        &B &94.31[0.40](0)&92.69[0.76](0.21)&95.00[1.12](26.44)\\\\
       &(100,101) &A &94.72[0.40](0)&93.40[0.55](0.19)&92.48[0.68](26.08)\\
       &        &B &94.21[0.40](0)&93.46[0.55](0.21)&92.62[0.68](26.44)\\\\
       &(199,200)&A &94.46[0.40](0)&93.51[0.45](0.19)&93.06[0.48](26.08)\\

       &       & B &93.92[0.40](0)&93.46[0.45](0.21)&92.71[0.48](26.44)& \\
       \hline

 \multicolumn{7}{c}{$\epsilon =\log(n)/n^{1/2} $}\\
\hline
$100$ & $(1,2)$  & A &$79.34[0.58](0.24)$&$72.51[1.05](88.03)$&$100.00[1.44](99.97)$\\
       &        &B &$78.94[0.58](0.24)$&$75.06[1.05](87.29)$&$66.67[1.27](99.97)$\\\\
       &$(50,51)$ &A &$78.51[0.58](0.24)$&$73.68[0.81](88.03)$&$100.00[0.87](99.97)$\\
       &        & B &$79.20[0.58](0.24)$&$72.46[0.80](87.29)$&$100.00[0.82](99.97)$\\\\
       &$(99,100)$&A &$78.81[0.58](0.24)$&$75.86[0.65](88.03)$&$100.00[0.70](99.97)$\\
       &       & B &$78.73[0.58](0.24)$&$75.30[0.66](87.29)$&$100.00[0.67](99.97)$\\\\
$200$ & $(1,2)$  & A &$82.64[0.41](0)$ &$69.26[0.80](70.75)$&$100.00[1.23](99.96)$\\
       &        & B &$82.51[0.41](0)$&$70.91[0.79](71.16)$&$100.00[0.83](99.96)$\\\\
       &$(50,51)$ &A &$83.03[0.41](0)$&$74.50[0.57](70.75)$&$50.00[0.94](99.96)$\\
       &        & B &$82.74[0.41](0)$&$75.24[0.57](71.16)$&$50.00[0.64](99.96)$\\\\
       &$(99,100)$&A &$82.66[0.41](0)$&$77.88[0.45](70.75)$&$100.00[0.50](99.96)$\\
       &       & B &$82.54[0.41](0)$&$79.96[0.45](71.16)$&$50.00[0.47](99.96)$\\
\hline
\end{tabular}
\label{table2}
\end{table}
\clearpage
\newpage

\subsection{Real Data Example}

We use the  affiliation network dataset in \cite{Sundaresan:Fischhoff:  Dushoff:  Rubenstein : 2007} as a data example. As discussed in \cite{Granovetter(2017)}, it remains an interesting issue that the animals should have some sort of privacy rights. In some ways, society has already begun moving in that direction. This network dataset  is based on a study of a community of $28$ Grevy's zebras.  \cite{Sundaresan:Fischhoff:  Dushoff:  Rubenstein : 2007} showed that Grevy's zebra individuals are more selective in their choices of associates, tending to form bonds with others in the same reproductive state.
In the dataset, Grevy's zebras are labelled from $1$ to $28$, and $111$ edges with finite weight $q=3$. The edge weight of $0$ denotes that a pair of zebras never appeared during the study,
the edge weight of $1$ denotes that a pair of zebras appeared together at least once, while the edge weight of $2$ indicates a statistically significant tendency of pairs to appear together.
On the other hand, the estimate $\bm{\hat{\alpha}}$ does not exist when the degree of a vertex is zero. Hence we removed the vertex $8$ whose degree is zero before analysis.
The network with the left $27$ vertices is shown in Figure \ref{g1(1)}.
We chose the privacy parameter $\varepsilon$ as $1$. Figure \ref{r22} reports the scatter  plots of noisy degree sequence $\bar{\textbf{d}}$ vs the estimates $\bm{\hat{\alpha}}$  for the $27$ Grevy's zebra dataset.
From Figure \ref{r22}, the value of $\bm{\hat{\alpha}}$ increases as the number of $\bar{\textbf{d}}$ increases.
Furthermore, the estimates can reveal a  trend in these zebras' choices of associates. The larger the estimates $\bm{\hat{\alpha}}$, the more zebras have associates or  the higher the frequency of pairs to appear together. As shown in Figure \ref{r22}, the number of zebras' associates and the frequency of pairs to appear together are more and higher under the case that $\hat{\bm \alpha}$ is around zero.
Table \ref{tabr}  reports the estimates, the $95\%$ confidence interval, the corresponding standard errors and the noisy degree sequence. In Table \ref{tabr}, the larger estimates correspond to the larger noisy degrees.
The largest degree is $23$ for vertex $4$, which also has the largest  estimate $0.447$. On the other hand, the  vertex for $22$ with the smallest estimate $-2.260$, has degree $2$ in Table \ref{tabr}.
\begin{figure}[!htb]
\centering
\label{g1(1)}
{
\setlength{\abovecaptionskip}{-3cm}
\includegraphics[height=0.7\textwidth=0.7]{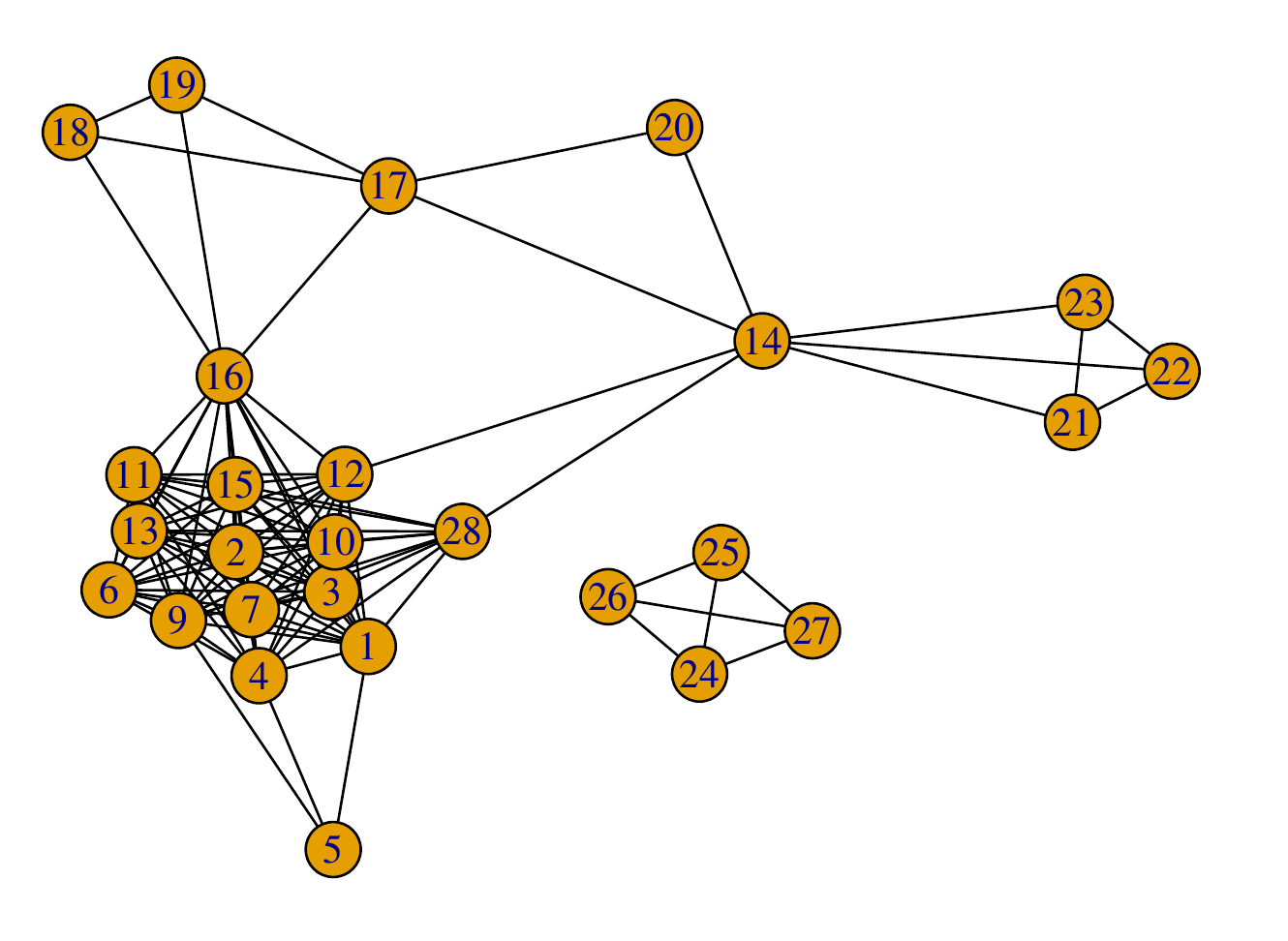}
}
\caption{the Sundaresan, Fischoff, Dushoff, Rubenstein Zebra Affiliation network becomes a community network of $27$ Greavey's zebras with $111$ edges, after removing the isolated vertex $8$. There are $3$ edge weights representing the trend of a pair of zebras appeared together during the study: $0$ means none; $1$ means at least once; $2$ means very significantly.}
\label{g1(1)}
\end{figure}
\begin{figure}[!htb]
\centering
\label{r22}
 \subfigure[$n=27,\epsilon=1$.]{
\setlength{\abovecaptionskip}{0.cm}
\setlength{\belowcaptionskip}{-0.cm}
\includegraphics[height=5in, width=5in]{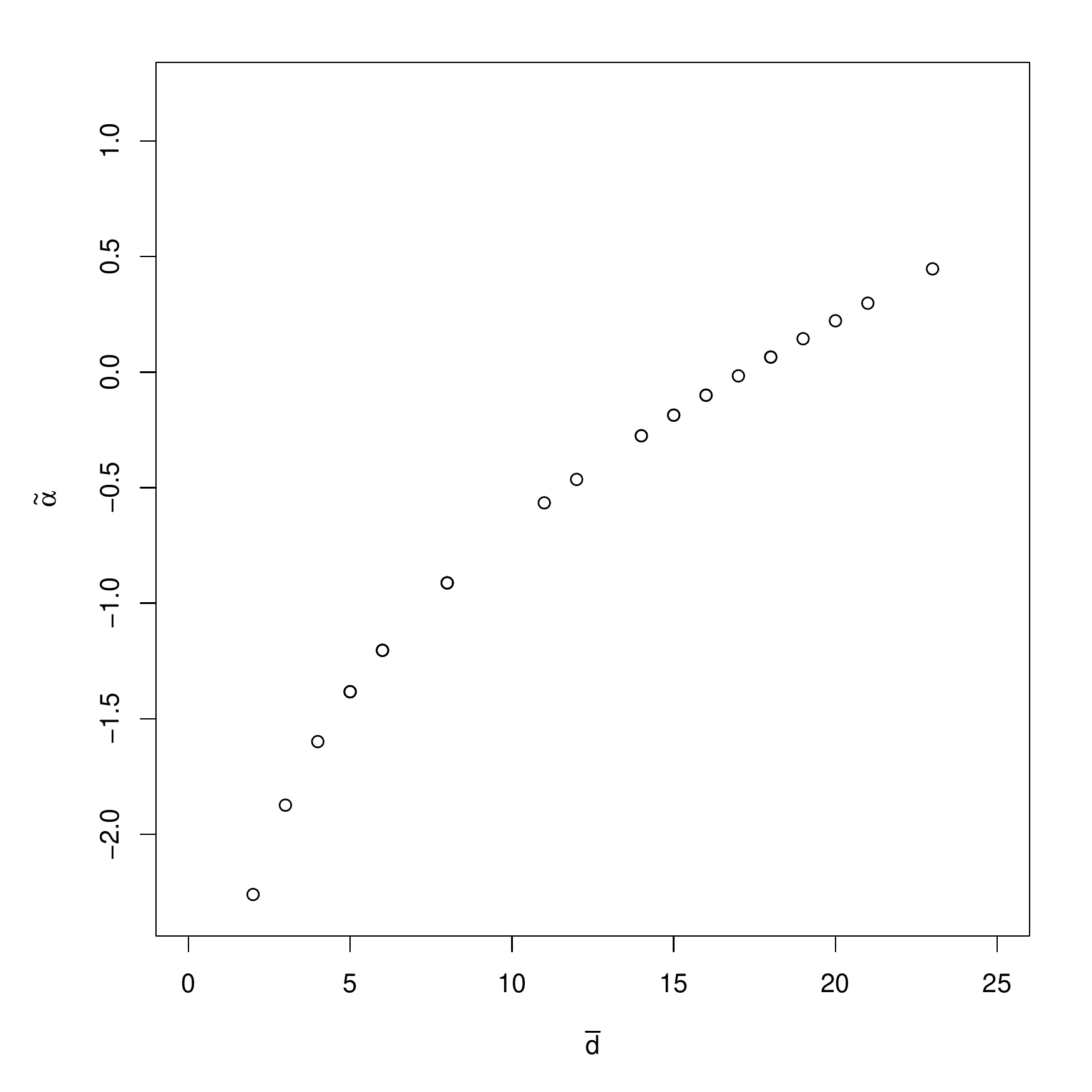}
}
\caption{The scatter  plots ($\bar{\textbf{d}}$ vs $\bm{\hat{\alpha}}$). The $\bar{\textbf{d}}$ denotes the noisy degree sequence and $\bm{\hat{\alpha}}$ denotes the corresponding estimate. }
\label{r22}
\end{figure}

\renewcommand\arraystretch{1}
\begin{table}[!ht]
 \centering
\caption{Sundaresan, Fischoff, Dushoff, Rubenstein Zebra Affiliation network dataset: the estimate $\bm{\hat{\alpha}}$, $95\%$ confidence intervals (in square brackets), and their standard errors (in parentheses).}

\scalebox{0.8}{
\begin{tabular*}{18.5cm}{llllll}
\hline
Vertex &$\hat{\alpha}_{i}$ &Degree &Vertex & $\hat{\alpha}_{i}$ & Degree \\
\hline
&&$\epsilon=1$\\
\hline
$1$ &  $ 0.065 [-0.477 ,  0.606] (0.276) $&$18$  &$16$ & $-0.186[-0.751, 0.379] (0.288)$ &$15$\\

$2$   &$0.298 [-0.229 ,  0.825] ( 0.269)$ &$21$      &$17$& $-1.204 [-1.994, -0.414] (0.403)$  &$6$\\
$3$   & $-0.276 [-0.851 ,   0.300] ( 0.294)$ &$14 $     &$18$ &$ -1.383 [-2.242 ,  -0.524] (0.438)$  &$5$\\
$4$  & $0.447 [-0.075,  0.968] (0.266)$ &$23 $     &$19 $&$-1.599 [-2.554 , -0.643] (0.488)$  &$4$\\
$5$  &$-0.912   [ -1.611 ,  -0.213] ( 0.356) $ &$8 $     &$20 $&$-1.204  [ -1.994, -0.414 ] (0.403) $ &$6$\\
$6 $  &$-0.186   [  -0.751,  0.379 ]  (0.288) $& $15$       &$21 $& $-1.383 [-2.242 , -0.524] (0.438)$  &$5$\\
$7 $ &$ -0.276 [ -0.851, 0.300]  (0.294) $&$14 $     &$22 $& $ -2.260    [ -3.611 ,  -0.910 ] (0.689)$  &$2$\\
$9$  &$0.065 [ -0.477, 0.606] ( 0.276) $&$18 $       &$23$ &$ -0.912   [ -1.611, -0.213 ] (0.356) $ &$8$\\
$10 $ & $0.144 [ -0.391 ,  0.680] (0.273)$ &$ 19 $     &$24$ &$ -1.874  [-2.975,  -0.773 ] (0.562)$ &$3$\\
$11$ &$-0.100  [-0.656, 0.456] (0.284)$ &$ 16 $     &$25$ &$-0.464  [ -1.067, 0.138] (0.307)$  &$12$\\
$12$  &$-0.016 [-0.565 ,   0.532] (0.280)$ &$17$      &$26$ & $-1.204  [ -1.994,  -0.414 ] (0.403)$ &$6$\\
$13$ &$ -0.100 [ -0.656 ,   0.456] (0.284)$ &$16 $    &$27$ &$ -0.912 [-1.611, -0.213] (0.356) $ &$8$\\
$14 $&$-1.383 [ -2.242, -0.524] (0.438) $ &$5 $    &$28$ &$-0.566  [-1.186  ,  0.054] (0.316)$ &$11$\\
$15$ & $0.222 [ -0.309, 0.753] (0.271)$ &$20$\\
\hline
\end{tabular*}}
\label{tabr}
\end{table}

\newpage
\section{Summary and Future Study}\label{s4}
\par In this paper, we have established the uniform consistency and asymptotic normality of the Z-estimator in the generalized $\beta$-model involving noisy degree sequence $\bar{\textbf{d}}=\textbf{d}+\textbf{e}$, where $\textbf{d}$ is the sufficient statistic and $ \textbf{e}$  are some noises from discrete Laplace distribution. By using the Newton-Kantorovich theorem, we try to ignore adding noisy process, and obtain the existence and consistency of the $Z$-estimator satisfying the equation (\ref{eq1}).  Furthermore, we give some simulation results to illustrate the authenticity of our obtained results under a non-denoised process. Although the simulation results in Section \ref{s3} show that the approximation of asymptotic normality behaves well under  certain conditions, our theoretical results may be improved but not only in the generalized $\beta$-model.

It should be noted that the discrete Laplace random variable is the difference of two i.i.d. geometric distributed random variables, see Proposition $3.1$ in \cite{Inusah2006}. The geometric distribution as the class of infinitely divisible distribution is a special case of discrete compound Poisson distribution. The difference of geometric noise-addition mechanism can be flexibly extended to the difference between two i.i.d. (or independent) discrete compound Poisson random variables, see Definition $4.2$ of \cite{Zhang14}. In fact, the difference of two independent discrete compound Poisson random variables follows the infinitely divisible distributions with integer support, see Chapter IV of \cite{Steutel2003}. The frequently employed discrete Laplace noise $\{e_{i}\}$ in differential privacy, may be further optimally selected from other flexible discrete distributions to achieve effective privacy protection in the future.

\section{Appendix: Proofs}\label{s5}

\subsection{ Proof of Theorem \ref{th21}}
The proof of Theorem \ref{th21} is based on two steps.

\textbf{Step 1}； we need two lemmas below.

\begin{lemma}\label{lem1}
{\rm Consider the discrete Laplace mechanism with ${\lambda _n} = \exp ( -\epsilon_n/2)$, if $\varepsilon_{n}=\Omega(\sqrt{\log n})$ and $e^{Q_{n}}=o((n/\log n)^{\frac{1}{12}})$, then as $n$ goes to infinity, for any fixed $r\geq1,$
\[(\frac{{\overline {{d_1}}  - E({d_1})}}{\sqrt{{v_{11}}}}, \cdots ,\frac{{\overline {{d_r}}  - E({d_r})}}{\sqrt{{v_{rr}}}})\stackrel {d} { \rightarrow } {N_r}(0,\mathbf{I}_r).\]}
\end{lemma}
Lemma \ref{lem1} indicates that the components of $(\bar{d_{1}}-E(d_{1}),\cdots,\bar{d_{r}}-E(d_{r}) )$ are asymptotically independent and normally distributed with variances $v_{11},\cdots,v_{rr}, $ respectively.

\begin{lemma}\label{lem2}
 {\rm(\cite{Karwa2016}, Proposition E) Let $e_{1},\cdots, e_{n}$ be i.i.d random variables drawn from
discrete Laplace distribution with  probability mass function defined by
$$P(e_{i}=e)=\frac{1-\lambda}{1+\lambda}\lambda^{|e|}, e\in \mathcal{Z}, \lambda\in(0,1).$$
Then we have $E(e_{i})=0$ and $\text{Var}(e_{i})=\cfrac{2\lambda}{(1-\lambda)^{2}}.$ Moreover,
\[
    \begin{split}
    P(|e_{i}|>c)=\frac{2\lambda^{[c]+1}}{1+\lambda},~~P(\max_{i}&|e_{i}|>c)=1-\left(1-\cfrac{2\lambda^{[c]+1}}{1+\lambda}\right)^{n}~~ \text{and}~~E|e_{i}|=\frac{2\lambda}{1-\lambda^{2}}.
\end{split}
    \]
}
 \end{lemma}

\noindent{{\bf Proof of Lemma \ref{lem1}}
\par By $\bar{d}_{i}=d_{i}+e_{i}$, we can analysis the asymptotic normality of the following proposition in two parts, i.e., $(d_{i}-E(d_{i}))/v^{1/2}_{ii}$ and $(e_{i}-E(e_{i}))/v^{1/2}_{ii}$. On the one hand, \cite{Yan:Zhao:Qin:2015} have verified the result of the first part by Liapounov's central limit theorem [\cite{Chung2001}]. On the other hand, we can easily obtain the stochastic order of the second part by Chebyshev inequality.

Let $\bar{d}_{i}=d_{i}+e_{i},\ \ i=1,\cdots, r,$ then
\begin{equation}\label{eq51}
\frac{\bar{d}_{i}-E(d_{i})}{\sqrt{v_{ii}}}=\frac{d_{i}-E(d_{i})}{\sqrt{v_{ii}}}+\frac{e_{i}}{\sqrt{v_{ii}}}, \ \ i=1,\cdots, r .
\end{equation}
%Here, since \cite{Yan:Qin:Wang:2016b}(see p.277) have proved that  $ \cfrac{d_{i}-E(d_{i})}{\sqrt{v_{ii}}},i=1,\cdots, r$ are  asymptotically standard normal as $n\rightarrow\infty,$
Now, we only discuss the property of $\cfrac{e_{i}}{\sqrt{v_{ii}}}.$ In fact, by Chebyshev inequality and Lemma \ref{lem2}, for any constant $a>0,$ as $n$ goes to infinity, we have
\[
    \begin{split}
P(\big|\frac{e_{i}}{\sqrt{v_{ii}}}\big|\geq a)&=P(\big|e_{i}\big|\geq a\sqrt{v_{ii}})\leq\cfrac{\text{Var}(e_{i})}{a^{2}v_{ii}}\\
[\text{by}~(\ref{eq2})]~~&\leq\cfrac{2(1+e^{Q_{n}})}{a^{2}(n-1)}\times \cfrac{2e^{-\frac{\varepsilon_n}{2}}}{(1-e^{-\frac{\varepsilon_n}{2}})^{2}}\\
&\leq O\left(\frac{e^{Q_{n}}}{n}\right)\rightarrow 0.
 \end{split}
    \]}
by noticing that ${e^{ - \frac{{{\varepsilon _n}}}{2}}} < 1$ for all $n$.

\textbf{Step 2}: we apply the Newton-Kantorovich theorem (\cite{Gragg:Tapia:1974}) to obtain the existence and consistency of  the  estimator satisfying the equation (\ref{eq1}).

For a subset $C\subset R^{n},$ let $C^{0}$ and $\overline{C}$ denote the interior and closure of $C$ in $R^{n}$, respectively. Let $\Omega(\textbf{x},r)$ denote the open ball $\{\textbf{y}:\Vert \textbf{y}-\textbf{x}   \Vert<r\}$, and $\overline{\Omega(\textbf{x},r)}$ be its closure.
\begin{proposition}
{\rm(\cite{Gragg:Tapia:1974}) Let $F(\textbf{x})=(F_{1}(\textbf{x}),\cdots,F_{n}(\textbf{x}))^{T}$ be a function vector on $\textbf{x}\in R^{n}.$ Assume that the Jacobian matrix $F'(\textbf{x})$ is Lipschitz  continuous on an open convex set $D$ with the Lipschitz constant $\kappa.$ Given $\textbf{x}_{0}\in D,$ assume that $[F'(\textbf{x}_{0})]^{-1}$ exists,
\[
    \begin{split}
    &\Vert[F'(\textbf{x}_{0})]^{-1} \Vert_{\infty}\leq\aleph, \ \ \Vert[F'(\textbf{x}_{0})]^{-1}F'(\textbf{x}_{0})] \Vert_{\infty}\leq\delta, \ \ h=2\aleph\kappa\delta\leq1,\\
    &\Omega(\textbf{x}_{0}, t^{*})\subset D^{0}, t^{*} :=\frac{2}{h}(1-\sqrt{1-h})\delta=\frac{2}{1+\sqrt{1-h}}\delta\leq2\delta,
    \end{split}
    \]
where $\aleph$ and $\delta$ are positive constants that may depend on $\textbf{x}_{0}$ and the dimension $n$ of $\textbf{x}_{0}.$ Then the Newton iteration $\textbf{x}_{k+1}=\textbf{x}_{k}-[F'(\textbf{x}_{k})]^{-1}F(\textbf{x}_{k})$ exists and $\textbf{x}_{k}\in \Omega(\textbf{x}_{0},t^{*})\subset D^{0}$ for all $k\geq 0; \hat{\textbf{x}}=\lim \textbf{x}_{k}$ exists, $\hat{\textbf{x}}\in\overline{\Omega(\textbf{x}_{0},t^{*})}\subset D$ and $F(\hat{\textbf{x}})=0.$

}
\end{proposition}

 Besides,  we also need the following four lemmas to prove Theorem \ref{th21}. Specifically,  Lemmas \ref{lem3}--\ref{lem5} are served for Newton-Kantorovich theorem, and Lemma \ref{lem6} is based on Hoeffding's inequality and Lemma \ref{lem2}.
 \begin{lemma}\label{lem3}
{\rm(\cite{Yan:Xu:2013}) If $V_{n}\in \mathcal{L}_{n}(m,M)$, and $n$ is large enough, then
$$\Vert V_{n}^{-1}-\bar{S}_{n} \Vert\leq\frac{cM^{2}}{m^{3}(n-1)^{2}},$$
where $(\bar{S}_{n})_{ij}:=\cfrac{\delta_{ij}}{v_{ii}}-\cfrac{1}{v_{..}},v_{..}:=\sum^{n}_{i=1}v_{ii}, c$ is a constant that not depends on $M, m$ and $n.$
}
\end{lemma}
\par To establish the form $\Vert [F'(\bm{\alpha})]^{-1}F(\bm{\alpha})\Vert_{\infty}\leq\delta$ in Theorem \ref{th21}, we first use a simple matrix $S_{n}=(s_{ij})$ to approximate $V^{-1}_{n}$. The upper bound of the approximation error is given below.
\begin{lemma}\label{lem4}
{\rm Assume that $\bm{\alpha}\in D$, where $D=\{\bm{\alpha}\in R^{n}: -Q_{n}\leq \alpha_{i}+\alpha_{j}\leq Q_{n}, \text{for}\ 1\leq i<j\leq n \}$. If $n$ is large enough, then
$$\Vert V_{n}^{-1}-{S}_{n} \Vert\leq\frac{c_{1}(1+e^{Q_{n}})^{3}}{(n-1)^{2}},$$
where $V_{n}:=F'(\bm{\alpha}),$ $(S_{n})_{ij}:=\cfrac{\delta_{ij}}{v_{ii}}, c_{1}$ is a constant that not depends on $e^{Q_{n}}$ and $n.$\\
{\bf Proof.} By (\ref{eq2}), $v_{..}=\sum^{n}_{i=1}v_{ii}\ge n(n-1) m$ and Lemma \ref{lem3}, we can easily obtain
\[
    \begin{split}
\Vert V^{-1}_{n}-S_{n}\Vert&\leq \Vert V^{-1}_{n}-\bar{S}_{n}\Vert +\Vert\bar{S}_{n}-S_{n}\Vert\\
&\leq \frac{cM^{2}}{m^{3}(n-1)^{2}}+\frac{1}{n(n-1)m}\leq(c+\frac{m^{2}}{M^{2}})\frac{M^{2}}{m^{3}(n-1)^{2}}\\
&\leq(c+1)\frac{M^{2}}{m^{3}(n-1)^{2}}=\frac{c_{1}(1+e^{Q_{n}})^{3}}{(n-1)^{2}},
     \end{split}
    \]
where $c_{1}$ is a constant that not depends on $e^{Q_{n}}$ and $n.$
}
\end{lemma}
\par To confirm the value of $\aleph$ in Newton-Kantorovich theorem, we use triangle inequality and Lemma \ref{lem4} to obtain the upper bound of $V_{n}^{-1}$, as follow.
\begin{lemma}\label{lem5}
{\rm Assume that $\bm{\alpha}\in D$, where $D=\{\bm{\alpha}\in R^{n}: -Q_{n}\leq \alpha_{i}+\alpha_{j}\leq Q_{n}, \text{for}\ 1\leq i<j\leq n \}$. If $n$ is large enough, then
$$\Vert V_{n}^{-1} \Vert_{\infty}\leq\frac{c_{2}(1+e^{Q_{n}})^{3}}{n-1},$$
where $V_{n}:=F'(\bm{\alpha}),$ $(S_{n})_{ij}:=\cfrac{\delta_{ij}}{v_{ii}}, c_{2}$ is a constant that not depends on $e^{Q_{n}}$ and $n.$\\
{\bf Proof.} By (\ref{eq2}) and Lemma \ref{lem4}, we  obtain
\[
    \begin{split}
\Vert V_{n}^{-1} \Vert_{\infty}&\leq\Vert V^{-1}_{n}-S_{n}\Vert_{\infty}+\Vert S_{n}\Vert_{\infty}
\leq \frac{c_{1}n(1+e^{Q_{n}})^{3}}{(n-1)^{2}}+\frac{2(1+e^{Q_{n}})}{n-1}
\leq\frac{c_{2}(1+e^{Q_{n}})^{3}}{n-1},
     \end{split}
    \]
where $c_{2}$ is a constant that not depends on $e^{Q_{n}}$ and $n.$
}
\end{lemma}
\par The following lemma  guarantees that the upper bound of $\Vert\bar{\textbf{d}}-E(\textbf{d})\Vert_{\infty}$ is the magnitude of $(n\log n)^{1/2}$.
\begin{lemma}\label{lem6}
{\rm Let $\kappa_{n}=2(q-1)\sqrt{(n-1)\log (n-1)}.$ If $\varepsilon_{n}=\Omega(\sqrt{\log n}),$ then with probability approaching one as $n\rightarrow\infty$,
\begin{equation}\label{eq52}
\max_{1\leq i\leq n}\vert \bar{d}_{i}-E(d_{i})\vert \leq 2(q-1)\sqrt{(n-1)\log(n-1)}.
    \end{equation}
{\bf Proof.} Let $\kappa_{n}=2(q-1)\sqrt{(n-1)\log (n-1)},$ then
\begin{equation}\label{eq53}P(\max_{i}\mid\bar{d}_{i}-E(d_{i})\mid\geq\kappa_{n})\leq P(\max_{i}\mid d_{i}-E(d_{i})\mid\geq\frac{\kappa_{n}}{2})+ P(\max_{i}\mid e_{i}\mid\geq\frac{\kappa_{n}}{2}).\end{equation}
Here, the inequality (\ref{eq53}) is divided into two parts and is given respective discussions. For the first part, we have
$$P(\max_{i}\mid d_{i}-E(d_{i})\mid\geq\cfrac{\kappa_{n}}{2})\leq \sum_{i}P(\mid d_{i}-E(d_{i})\mid\geq\cfrac{\kappa_{n}}{2}).$$
By Hoeffding's inequality, it implies
\begin{equation}\label{eq54}
    \begin{split}
    P(\mid d_{i}-E(d_{i})\mid\geq\cfrac{\kappa_{n}}{2})&\leq 2\exp\left(-\frac{2\left(\frac{\kappa_{n}}{2}\right)^{2}}{(n-1)(q-1)^{2}}\right) \\
    &=2\exp\left(-\frac{2(q-1)^{2}(n-1)\log(n-1)}{(n-1)(q-1)^{2}}\right)\\&=\frac{2}{(n-1)^{2}}.
\end{split}
   \end{equation}
 For the second part, by definition of maximum and ${\lambda _n} = \exp ( -\epsilon_n/2 )$,  we obtain
\begin{align*}
    P(\max_{i}\mid e_{i}\mid\geq\cfrac{\kappa_{n}}{2})&=1-\prod^{n}_{i=1}P(\mid e_{i}\mid\leq\cfrac{\kappa_{n}}{2})\\
(\text{By Lemma~\ref{lem2}})~&=1-(1-\frac{2{\lambda_n}^{[\frac{\kappa_{n}}{2}]+1}}{1+{\lambda_n}})^{n}\\
    &=1-(1-\frac{2e^{-\epsilon_{n}([\frac{\kappa_{n}}{2}]+1)/2}}{1+e^{-\frac{\epsilon_{n}}{2}}})^{n}.
\end{align*}

 For $x\in(0,1),$ $f(x)=1-(1-x)^{n}$ is an increasing function on $x$. Then we get
 $$1-(1-\frac{2e^{-\epsilon_{n}([\frac{\kappa_{n}}{2}]+1)/2}}{1+e^{-\frac{\epsilon_{n}}{2}}})^{n}\leq 1-(1-2e^{-\epsilon_{n}([\frac{\kappa_{n}}{2}]+1)/2})^{n}.$$
 For $x\in (0,1),$  $(1-x)^{n}\geq1-nx.$  Then we also get
  $$1-(1-\frac{2e^{-\epsilon_{n}([\frac{\kappa_{n}}{2}]+1)/2}}{1+e^{-\frac{\epsilon_{n}}{2}}})^{n}\leq 1-(1-2ne^{-\epsilon_{n}([\frac{\kappa_{n}}{2}]+1)/2})=2ne^{-\epsilon_{n}([\frac{\kappa_{n}}{2}]+1)/2}\leq2ne^{-\epsilon_{n}[\frac{\kappa_{n}}{2}]/2}.$$
 By $\epsilon_{n}=\Omega(\sqrt{\log n}),$ we obtain $\epsilon_{n}\geq c\sqrt{\log n} \geq c\sqrt{\frac{\log n}{n}}$, and thus
  $\epsilon_{n}[\frac{\kappa_{n}}{2}]\geq 4\log n.$  This implies
\begin{equation}\label{eq55}P(\max_{i}|e_{i}|\geq \cfrac{\kappa_{n}}{2})\leq \cfrac{2n}{n^{2}}.\end{equation}
Therefore, by (\ref{eq54}) and (\ref{eq55}), with probability approaching one as $n\rightarrow\infty$, we have
$$\max_{i}|\bar{d_{i}}-E(d_{i})|\leq 2(q-1)\sqrt{(n-1)\log(n-1)}.$$
}
\end{lemma}

\noindent{{\bf Proof of Theorem \ref{th21}} In the Newton's iterative step, putting the initial value $\bm{\alpha}^{0}:=\bm{\alpha}.$ Let $V_{n}=F'(\bm{\alpha})\in \mathcal{L}_{n}(m,M)$ and $W_{n}=V_{n}^{-1}-S_{n}.$ Let $F(\bm{\alpha})=\bar{\textbf{d}}-E(\textbf{d})$, by Lemma \ref{lem4} and (\ref{eq52}), we get
\[
    \begin{split}
\Vert [F'(\bm{\alpha})]^{-1}F(\bm{\alpha})\Vert_{\infty}&\leq n\Vert W_{n}\Vert\Vert F(\bm{\alpha})\Vert_{\infty}+\max_{i}\frac{\vert F_{i}(\bm{\alpha})\vert}{v_{ii}}\leq (n\Vert W_{n} \Vert+\frac{1}{v_{ii}})\Vert F(\bm{\alpha}) \Vert_{\infty}\\
&\leq \left(\frac{c_{1}(1+e^{Q_{n}})^{3}}{(n-1)^{2}}+\frac{2(1+e^{Q_{n}})}{n-1}\right)\Vert F(\bm{\alpha}) \Vert_{\infty}\\
&\leq \left(\frac{c_{1}(1+e^{Q_{n}})^{3}}{(n-1)^{2}}+\frac{2(1+e^{Q_{n}})}{n-1}\right)\times 2(q-1)\sqrt{(n-1)\log(n-1)}\\
&\leq c_{3}(1+e^{Q_{n}})^{3}\sqrt{\frac{\log (n-1)}{n-1}},
\end{split}
    \]
where $c_{3}$ is a constant.

Combining Lemma \ref{lem5},  we can set $\delta= c_{3}(1+e^{Q_{n}})^{3}\sqrt{\frac{\log (n-1)}{n-1}}$ and $\aleph=\frac{c_{2}(1+e^{Q_{n}})^{3}}{n-1}$ in Newton-Kantorovich theorem. \par
Next, we indicate that the Jacobian matrix $F'(\bm{\alpha})$ is Lipschitz  continuous with $\kappa=4(q-1)^{3}(n-1)$. Here, our method is similar to \cite{Yan:Qin:Wang:2016b}.\par
Let $$\textmd{g}_{ij}(\bm \alpha)=(\frac{\partial^{2}F_{i}}{\partial\alpha_{1}\partial\alpha_{j}},\cdots,\frac{\partial^{2}F_{i}}{\partial\alpha_{n}\partial\alpha_{j}} )^{T}.$$
By some computations, we have
\[
    \begin{split}
&\frac{\partial ^{2}F_{i}}{\partial\alpha_{i}^{2}}=\sum_{j=1;j\neq i}^{n}\frac{(1/2)\sum_{k\neq l,a}(k-l)^2(k+l-2a)e^{(k+l+a)(\alpha_{i}+\alpha_{j})}}{(\sum_{a=0}^{q-1}e^{a(\alpha_{i}+\alpha_{j})})^{3}},\\
&\frac{\partial ^{2}F_{i}}{\partial\alpha_{j}\alpha_{i}}=\frac{(1/2)\sum_{k\neq l,a}(k-l)^2(k+l-2a)e^{(k+l+a)(\alpha_{i}+\alpha_{j})}}{(\sum_{a=0}^{q-1}e^{a(\alpha_{i}+\alpha_{j})})^{3}}.
\end{split}
    \]
As $\sum_{k\neq l,a}e^{(k+l+a)(\alpha_{i}+\alpha_{j})}\leq(\sum_{a=0}^{q-1}e^{a(\alpha_{i}+\alpha_{j})})^{3},$ uniformly we have
\begin{equation}\label{eq56}
\vert\cfrac{\partial ^{2}F_{i}}{\partial\alpha_{i}^{2}}\vert\leq(q-1)^{3}(n-1),\ \ \vert\frac{\partial ^{2}F_{i}}{\partial\alpha_{j}\alpha_{i}}\vert\leq (q-1)^{3}.
    \end{equation}
Thus $\Vert\textmd{g}_{ii}(\bm\alpha)\Vert_{1}\leq 2(n-1)(q-1)^{3}.$

If $i\neq j$ and $k\neq i,j,$
 $$\frac{\partial^{2}F_{i}}{\partial\alpha_{k}\partial\alpha_{j}}=0.$$
Then we get $\Vert\textmd{g}_{ij}(\bm\alpha)\Vert_{1}\leq 2(q-1)^{3},i\neq j.$  Therefore, by the mean-value for vector-valued functions (\cite{Lang1993}, p.341),  for a vector $\textbf{v},$
\begin{align*}
\max_{i}\left\{\sum_{j=1}^{n}\left[\frac{\partial F_{i}}{\partial \alpha_{j}}(\textbf{x})-\frac{\partial F_{i}}{\partial \alpha_{j}}(\textbf{y})\right]v_{j}\right\}&\leq \Vert   \textbf{v}\Vert_{\infty}\max_{i}\sum_{j=1}^{n}\left|\frac{\partial F_{i}}{\partial \alpha_{j}}(\textbf{x})-\frac{\partial F_{i}}{\partial \alpha_{j}}(\textbf{y})\right|\\
&=\Vert   \textbf{v}\Vert_{\infty}\max_{i}\sum_{j=1}^{n}\left|\int^{1}_{0}\textmd{g}_{ij}(t  \textbf{x}+(1-t)\textbf{y})(\textbf{x}-\textbf{y})dt\right|\\
&\leq 4(q-1)^{3}(n-1)\Vert   \textbf{v}\Vert_{\infty}\Vert  \textbf{x}-\textbf{y}\Vert_{\infty}.
\end{align*}
Thereby, we have
\[
    \begin{split}
 h&=2\aleph\kappa\delta=8\times\cfrac{c_{2}(1+e^{Q_{n}})^{3}}{n-1}\times(q-1)^{3}(n-1)\times c_{3}(1+e^{Q_{n}})^{3}\sqrt{\frac{\log (n-1)}{n-1}},\\
&=8c_{2}c_{3}(q-1)^{3}(1+e^{Q_{n}})^{6}\sqrt{\frac{\log (n-1)}{n-1}}=O(e^{6Q_{n}}\sqrt{\frac{\log n}{n}})=o(1).
\end{split}
    \]
Thus, all conditions in the Newton-Kantorovich theorem are satisfied. Since the inequality (\ref{eq52}) holds with probability approaching one, then (\ref{eq3}) is  fulfilled.
}
\subsection{ Proof of Theorem \ref{th22}}
\par The aim of proving Theorem \ref{th22} is to establish the following equation
\[
    \begin{split}
(\bm{\hat{\alpha}}-\bm{\alpha})_{i}=[S_{n}(\bar{\textbf{d}}-E(\textbf{d}))]_{i}+o_{p}(n^{-1/2})
    .
\end{split}
    \]
This will follow directly the from $[W_{n}\{\bar{\textbf{d}}-E(\textbf{d})\}]_{i}=o_{p}(n^{-1/2}), n\rightarrow\infty,$ and by Theorem  \ref{th21}. To this end, we introduce the following lemma.
\begin{lemma}\label{lem7}
{\rm Assume that $\bm{\alpha}\in D$, where $D=\{\bm{\alpha}\in R^{n}: -Q_{n}\leq \alpha_{i}+\alpha_{j}\leq Q_{n}, \text{for}\ 1\leq i<j\leq n \}$.  Let $\varepsilon_{n}=\Omega(\sqrt{\log n}), e^{Q_{n}}=o((n/\log n)^{\frac{1}{12}})$, and $U_{n}=\text{cov}[W_{n}\{\bar{\textbf{d}}-E(\textbf{d})\}].$ Then
$$[W_{n}\{\bar{\textbf{d}}-E(\textbf{d})\}]_{i}=o_{p}(n^{-1/2}).$$
{\bf Proof.} Let $\overline{V}_{n}=\text{cov}\{\bar{\textbf{d}}-E(\textbf{d})\}, V_{n}=\text{cov}\{\textbf{d}-E(\textbf{d})\}$ and $E_{n}=\text{cov}(\textbf{e}).$ For $1\leq i\leq n,$ the random variables $d_{i}$ and $e_{i}$ are mutually independent, then
\[
    \begin{split}
\text{cov}(\bar{d}_{i}-E(d_{i}),\bar{d}_{j}-E(d_{j}))&=\text{cov}(d_{i}+e_{i}-E(d_{i}), d_{j}+e_{j}-E(d_{j}))\\
&=\text{cov}(d_{i}-E(d_{i}), d_{j}+e_{j}-E(d_{j}))+\text{cov}(e_{i}, d_{j}+e_{j}-E(d_{j}))\\
&=\text{cov}(d_{i}-E(d_{i}), d_{j}-E(d_{j}))+\text{cov}(d_{i}-E(d_{i}),e_{j})\\&+\text{cov}(e_{i}, d_{j}-E(d_{j}))
+\text{cov}(e_{i}, e_{j})\\
&=\text{cov}(d_{i}-E(d_{i}), d_{j}-E(d_{j}))+\text{cov}(e_{i}, e_{j}).
\end{split}
    \]
Two cases are discussed:\\
Case 1. If $i\neq j$, then $\text{cov}(\bar{d}_{i}-E(d_{i}),\bar{d}_{j}-E(d_{j}))=\text{cov}(d_{i}-E(d_{i}), d_{j}-E(d_{j}));$\\
Case 2. If $i=j$, then $\text{cov}(\bar{d}_{i}-E(d_{i}),\bar{d}_{j}-E(d_{j}))=\text{Var}(d_{i}-E(d_{i}))+\text{Var}(e_{i}).$\\
Thus the elements of the  matrix $\overline{V}_{n}$ are denoted by $\bar{v}_{ij}=v_{ij},\bar{v}_{ii}=v_{ii}+\text{Var}(e_{i}),\  1\leq i\neq j \leq n.$
Let $U_{n}=\text{cov}[W_{n}\{\bar{\textbf{d}}-E(\textbf{d})\}]$ with $W_{n}=V^{-1}_{n}-S_{n},$ then $$U_{n}=W_{n}\overline{V}_{n}W^{T}_{n}=W_{n}(V_{n}+E_{n})W_{n}^{T}=W_{n}V_{n}W_{n}^{T}+W_{n}E_{n}W_{n}^{T}.$$
On the one hand, $W_{n}V_{n}W^{T}_{n}=(V_{n}^{-1}-S_{n})-S_{n}(I_{n}-V_{n}S_{n}),$ where $n\times n$ matrix $I_{n}$ is an  identity matrix.\\
By (\ref{eq2}), we obtain
\begin{equation}\label{eq57}
    \begin{split}
\big| \{S_{n}(I_{n}-V_{n}S_{n})\}_{ij}\big|=\big|\cfrac{(\delta_{ij}-1)v_{ij}}{v_{ii}v_{jj}}\big|\leq\frac{2q^{2}(1+e^{Q_{n}})^{2}}{(n-1)^{2}}.
\end{split}
    \end{equation}
By Lemma \ref{lem4} and (\ref{eq57}), we have
\[
    \begin{split}
   \Vert W_{n}V_{n}W^{T}_{n}\Vert \leq\frac{c_{1}(1+e^{Q_{n}})^{3}}{(n-1)^{2}}+\frac{2q^{2}(1+e^{Q_{n}})^{2}}{(n-1)^{2}}\leq O\left(\cfrac{e^{3Q_{n}}}{n^{2}}\right).
\end{split}
    \]
On the other hand, \[
    \begin{split}
\big\| (W_{n}E_{n}W^{T}_{n})_{ij}\big\|&=\big\|\sum^{n=1}_{k=1}w_{ik}e_{k}w_{kj}\big\|\leq\max_{k}|e_{k}| \sum^{n}_{k=1}\big|w_{ik} \big|\big|w_{kj}\big|\\
&\leq n\max_{k}|e_{k}| \Vert W_{n}\Vert^{2}\leq \cfrac{2ne^{-\frac{\varepsilon}{2}}}{(1-e^{-\frac{\varepsilon}{2}})^{2}}\times\frac{c^{2}_{1}(1+e^{Q_{n}})^{6}}{(n-1)^{4}}\leq O\left(\frac{e^{6Q_{n}}}{n^{3}}\right).
\end{split}
    \]
Hence, $\|U_{n}\|\leq O\left(\cfrac{e^{3Q_{n}}}{n^{2}}\right).$ Furthermore, by Chebyshev inequality, for any constant $a>0$, we get
\[
    \begin{split}
P\left(\cfrac{[W_{n}\{\bar{\textbf{d}}-E(\textbf{d})\}]_{i}}{n^{-1/2}}\geq a\right)&\leq P\left([W_{n}\{\bar{\textbf{d}}-E(\textbf{d})\}]_{i}\geq an^{-1/2}\right)\\
&\leq\cfrac{n[\text{cov}\{W_{n}(\bar{\textbf{d}}-E(\textbf{d}))\}]_{i}}{a^{2}}\leq O\left(\cfrac{e^{3Q_{n}}}{n}\right).
\end{split}
    \]
Then while $e^{Q_{n}}=o((n/\log n)^{\frac{1}{12}})$, $P\left(\cfrac{[W_{n}\{\bar{\textbf{d}}-E(\textbf{d})\}]_{i}}{n^{-1/2}}\geq a\right)\rightarrow 0,$ as $n\rightarrow\infty.$ Therefore, $$[W_{n}\{\bar{\textbf{d}}-E(\textbf{d})\}]_{i}=o_{p}(n^{-1/2}).$$
}
\end{lemma}
\noindent{{\bf Proof of Theorem \ref{th22}.} Let $\hat{r}_{ij}=\hat{\alpha}_{i}+\hat{\alpha}_{j}-\alpha_{i}-\alpha_{j}$. Under the conditions in Theorem \ref{th21}, we have from the consistency property
\[
    \begin{split}
    \max_{i\neq j}|\hat{r}_{ij}|=O_p(e^{3Q_{n}}\sqrt{\frac{\log n}{n}}).
    \end{split}
    \]
Let $u(t)= \sum\limits_{a=0}^{q-1}\cfrac{ae^{at}}{\sum^{q-1}_{k=0}e^{kt}}.$
For $i=1,\cdots, n,$ by the Taylor's expansion, we get
\[
    \begin{split}
   \bar{d_{i}}-E(d_{i})&=\sum_{j\neq i}(u(\hat{\alpha}_{i}+\hat{\alpha}_{j})-u(\alpha_{i}+\alpha_{j}))\\
   &=\sum_{j\neq i}[u'(\alpha_{i}+\alpha_{j})((\hat{\alpha}_{i}+\hat{\alpha}_{j})-(\alpha_{i}+\alpha_{j}))]+h_{i},
    \end{split}
    \]
where $h_{i}=(1/2)\sum_{j\neq i}u''(\hat{r}_{ij})[((\hat{\alpha}_{i}+\hat{\alpha}_{j})-(\alpha_{i}+\alpha_{j}))]^{2}$ and $\hat{r}_{ij}=t_{ij}(\alpha_{i}+\alpha_{j})+(1-t_{ij})(\hat{\alpha}_{i}+\hat{\alpha}_{j}), t_{ij}\in (0,1).$ Writing the above expressions into a matrix, we have
$$\bar{\textbf{d}}-E(\textbf{d})=V_{n}(\hat{\boldsymbol{\alpha}}-\boldsymbol{\alpha})+\textbf{h},$$
thus
$$\hat{\boldsymbol{\alpha}}-\boldsymbol{\alpha}=V^{-1}_{n}(\bar{\textbf{d}}-E(\textbf{d}))+V^{-1}_{n}\textbf{h},$$
where $\textbf{h}=(h_{1},\cdots, h_{n})^{T}.$

By (\ref{eq56}), we know $|h_{i}|\leq \cfrac{1}{2}(n-1)(q-1)^{3}\hat{r}^2_{ij}.$  Therefore,
\[
    \begin{split}
   |(V^{-1}_{n}\textbf{h})_{i}|=|(S_{n}\textbf{h})_{i}|+|(W_{n}\textbf{h})_{i}|\leq \max_{i}\frac{|h_{i}|}{v_{ii}}+\Vert W_{n}\Vert \sum_{i}|h_{i}|\leq O\left(e^{9Q_{n}}\frac{\log n}{n}\right),
    \end{split}
    \]
 If $e^{Q_{n}}=o\left(\cfrac{n^{1/18}}{(\log n)^{1/9}}\right),$ then $(V^{-1}_{n}\textbf{h})_{i}=o(n^{-1/2}).$

 By Theorem \ref{th21} and Lemma \ref{lem7}, for $i=1,\cdots, r,$ we have
\[
    \begin{split}
    (\bm{\hat{\alpha}}-\bm{\alpha})_{i}=[S_{n}(\bar{\textbf{d}}-E(\textbf{d}))]_{i}+o_{p}(n^{-1/2})
    =\frac{\bar{d}_{i}-E(d_{i})}{v_{ii}}+o_{p}(n^{-1/2}).
    \end{split}
    \]
Hence, Theorem \ref{th22} follows directly from Lemma \ref{lem1}. Finally, we conclude the proof by multiplying $\sqrt{{v_{ii}}}$ to left and right of the last display.

\section*{Acknowledgements}
The authors thank the Editor, an associated editor and one referee for their valuable comments that greatly improve the manuscript. The authors also thank Yujing Gao for her helpful suggestions. Yan is partially supported by the National Natural Science Foundation of China (No. 11771171) and the Fundamental Research Funds for the Central Universities (No. CCNU17TS0005). Fan is supported by Fundamental Research Funds for the Central Universities (Innovative Funding Project) (2019CX
ZZ071).

\section*{Supplementary material}
\par Supplement to ``Asymptotic Theory for Differentially Private Generalized $\beta$-models with Parameters Increasing." The supplementary material contains a brief introduction of \emph{skew discrete Laplace mechanism} in Subsection 2.2, and the QQ plots of parameter estimates with $\epsilon =\log(n)/n^{1/2}$ in Subsection 3.1.

%\section*{Supplementary Materials}

\newpage

\appendix
\vspace*{0.15cm}
\begin{center}
{\Large\bf Supplementary material to ``Asymptotic Theory for Differentially Private Generalized $\beta$-models
with Parameters Increasing"}
\vskip 1.2\baselineskip
{\large Yifan Fan,$^{1}$ \ \ Huiming Zhang,$^{^2*}$ \ \ Ting Yan$^{1}$}
\vskip 0.2\baselineskip
{\sl Department of Statistics, Central China Normal University, Wuhan, 430079, China$^{1}$}
\vskip 0.1\baselineskip
\vskip 0.2\baselineskip
{\sl School of Mathematical Sciences and Center for Statistical Science, Peking University, Beijing, 100871, China$^{2}$}
\vskip 0.1\baselineskip

\end{center}

This is a supplementary material in \cite{AFan2020} that contains a brief introduction of \emph{skew discrete Laplace mechanism} in Subsection 2.2, and the QQ plots of parameter estimates with $\epsilon =\log(n)/n^{1/2}$ in Subsection 3.1.
\vskip 5 pt

\section{ Skew discrete Laplace mechanism}
\par  We give a brief introduction to the \emph{skew discrete Laplace mechanism}. When the positive noises and negative noises arising with difference probability law, the skew discrete Laplace
distribution [\cite{AKozubowski2006}] as a discretization of non-symmetric Laplace
distribution, which is useful in applications to communications, engineering, and finance and economics, see \cite{AKotz2012} and references therein. A random variable $Z$ has the \emph{skew discrete Laplace distribution} with parameters $\lambda, \mu \in(0,1),$ if
$$
P(Z=z)=\frac{(1-\lambda)(1-\mu)}{1-\lambda \mu} \left\{\begin{array}{ll}{\lambda^{z},} & {z=0,1,2,3, \ldots} \\ {\mu^{|z|},} & {z=0,-1,-2,-3, \ldots}\end{array}\right.
$$
By the way, we have the following result as a remark which is a skew extension of Lemma 1 of \cite{AKarwa2016}. The Proposition \ref{Apro1} indicates that if skewness of the distribution is large (i.e. $\frac{\lambda }{\mu }$ or $\frac{\mu }{\lambda }$  is small than 1.). Thus, with $\lambda  \wedge \mu <\lambda$~(or $\mu$), the un-symmetry distribution will lead to bigger  $\varepsilon$- edge differentially private. So bigger  $\varepsilon$ implies less privacy protection. This proposition suggests us that in order to protect the topological information in the graph $G$, we had better choose some symmetric noise distributions.
\begin{proposition}\label{Apro1}
Let $f:\mathcal{G}\rightarrow \mathcal{Z}^{k},$ and let $Z_{1},\ldots, Z_{k}$ be independent and identically distributed discrete Laplace random variables with probability mass function defined in above. Then the algorithm which outputs $f(G)+(Z_{1},\ldots, Z_{k})$ with inputs $G$ is $\varepsilon$-edge differentially private, where
$$\varepsilon  = - {\Delta _G}(f)\log (\lambda  \wedge \mu ).$$
{\bf Proof.}
Let $g$ and $g^{\prime}$ be two graphs that differ by an
edge. The output of the mechanism is $f(G)+Z \in \mathbb{Z}^{k}$ and hence it is enough
to consider the probability mass function. For any $s \in \mathbb{Z}^{k},$ consider
\begin{align*}
P[f(G)+Z=s | G=g] &=P[Z=s-f(g) | G=g] \\
 &=\prod_{i=1}^{k}P\left[Z_{i}=s_{i}-f_{i}(g) | G=g\right] \\
&=C(\lambda ,\mu )\prod\limits_{i = 1}^k {\left( {{\lambda ^{\left| {{s_i} - {f_i}(g)} \right|}}{{\rm{1}}_{{s_i} - {f_i}(g) \in {\mathbb{Z}^ + }}} + {\mu ^{\left| {{s_i} - {f_i}(g)} \right|}}{{\rm{1}}_{{s_i} - {f_i}(g) \in {\mathbb{Z}^ - }}}} \right)}
\end{align*}
WLOG, we let $t := \frac{\lambda }{\mu }{\rm{ < 1}}$, thus $1>t>\lambda<\mu$. By triangle inequality for each $i$ we have $\left|s_{i}-f_{i}(g)\right|-\left|s_{i}-f_{i}\left(g^{\prime}\right)\right| \geq-\left|f_{i}(g)-f_{i}\left(g^{\prime}\right)\right| .$ Thus by the above derivation we obtain:
\begin{align*}
&~~~~ \frac{P[f(G)+Z=s | G=g]}{P[f(G)+Z=s | G=g^{\prime}]} =\prod\limits_{i = 1}^k {\frac{{{\lambda ^{\left| {{s_i} - {f_i}(g)} \right|}}{{\rm{1}}_{{s_i} - {f_i}(g) \in {\mathbb{Z}^ + }}} + {\mu ^{\left| {{s_i} - {f_i}(g)} \right|}}{{\rm{1}}_{{s_i} - {f_i}(g) \in {\mathbb{Z}^ - }}}}}{{{\lambda ^{\left| {{s_i} - {f_i}(g')} \right|}}{{\rm{1}}_{{s_i} - {f_i}(g') \in {\mathbb{Z}^ + }}} + {\mu ^{\left| {{s_i} - {f_i}(g')} \right|}}{{\rm{1}}_{{s_i} - {f_i}(g') \in {\mathbb{Z}^ - }}}}}}\\
& \le \prod\limits_{i = 1}^k \left( {{\lambda ^{\left| {{s_i} - {f_i}(g)} \right| - \left| {{s_i} - {f_i}(g')} \right|}}{{\rm{1}}_{{s_i} - {f_i}(g),{s_i} - {f_i}(g') \in {\mathbb{Z}^ + }}} + {{\left( {\frac{\lambda }{\mu }} \right)}^{\left| {{s_i} - {f_i}(g)} \right| - \left| {{s_i} - {f_i}(g')} \right|}}{{\rm{1}}_{{s_i} - {f_i}(g) \in {\mathbb{Z}^ + },{s_i} - {f_i}(g') \in {\mathbb{Z}^ - }}}} \right. \\
& \left.+ {\mu ^{\left| {{s_i} - {f_i}(g)} \right| - \left| {{s_i} - {f_i}(g')} \right|}}{{\rm{1}}_{{s_i} - {f_i}(g),{s_i} - {f_i}(g') \in {\mathbb{Z}^ - }}} + {\left( {\frac{\mu }{\lambda }} \right)^{\left| {{s_i} - {f_i}(g)} \right| - \left| {{s_i} - {f_i}(g')} \right|}}{{\rm{1}}_{{s_i} - {f_i}(g) \in {\mathbb{Z}^ - },{s_i} - {f_i}(g') \in {\mathbb{Z}^ + }}}\right)\\
&\le \prod\limits_{i = 1}^k {\left( {{\lambda ^{ - \left| {{f_i}(g) - {f_i}(g')} \right|}}{{\rm{1}}_{{s_i} - {f_i}(g),{s_i} - {f_i}(g') \in {\mathbb{Z}^ + }}} + {{\left( {\frac{\lambda }{\mu }} \right)}^{ - \left| {{f_i}(g) - {f_i}(g')} \right|}}{{\rm{1}}_{{s_i} - {f_i}(g) \in {\mathbb{Z}^ + },{s_i} - {f_i}(g') \in {\mathbb{Z}^ - }}}} \right.} \\
& \left.+ {\mu ^{ - \left| {{f_i}(g) - {f_i}(g')} \right|}}{{\rm{1}}_{{s_i} - {f_i}(g),{s_i} - {f_i}(g') \in {\mathbb{Z}^ - }}} + {\left( {\frac{\mu }{\lambda }} \right)^{\left| {{f_i}(g) - {f_i}(g')} \right|}}{{\rm{1}}_{{s_i} - {f_i}(g) \in {\mathbb{Z}^ - },{s_i} - {f_i}(g') \in {\mathbb{Z}^ + }}}\right)\\
&\leq {\left( {\lambda  \wedge \mu  \wedge t} \right)^{ - \sum\limits_{i = 1}^k {|{f_i}(g) - {f_i}(g')|} }} \le {{\lambda }^{ -\bigtriangleup_G(f)}} = {{\mathop{\rm e}\nolimits} ^\varepsilon}.
\end{align*}
where the third last inequality is by using the fact that $\alpha^{x}$ is a decreasing function for $x$ satisfying $0<x<1$, and the second inequality stems from the fact that the four terms with index function only appear one time in the product operation.

For this case, we can put
\[\varepsilon  =  -\bigtriangleup_G(f)\log{\lambda }.\]

Similarly, for $t := \frac{\lambda }{\mu }{\rm{ >1}}$. We get $1> \lambda >\mu <t^{-1}$. Then
\[\frac{{P[f(G) + Z = s|G = g]}}{{P[f(G) + Z = s|G = {g^\prime }]}} \le {\left( {\lambda \wedge \mu \wedge {t^{ - 1}}} \right)^{ - {\Delta _G}(f)}}\le{{\mu }^{ -\bigtriangleup_G(f)}} = {{\rm{e}}^\varepsilon }\]
which implies $\varepsilon  =- \bigtriangleup_G(f)\log \mu $.

Combining the two cases above, we immediately get
\[\varepsilon  =  - {\Delta _G}(f)[\log (\lambda ){1_{\lambda  \le \mu }} + \log (\mu ){1_{\lambda  \ge \mu }}]{\rm{ = }} - {\Delta _G}(f)\log (\lambda  \wedge \mu ).\]
 \end{proposition}

\par The following Lemma \ref{Alemm} evaluates the tail distribution of skew discrete Laplace, by employing similar arguments as in the proof of Theorem 2.1 and 2.2 (\cite{AFan2020}), we will have the similar consistency and asymptotic normality results under the skew discrete Laplace mechanism ${\lambda _n} = \exp ( - \frac{{{\varepsilon _n}}}{2})$ and ${\lambda _n}=O({\mu _n})$.
\begin{lemma}\label{Alemm}
 Let ${\tilde e}_{1},\cdots, {\tilde e}_{n}$ be i.i.d random variables drawn from
skew discrete Laplace distribution with probability mass function defined by
$$
P(Z={\tilde e}_1)=\frac{(1-\lambda)(1-\mu)}{1-\lambda \mu} \left\{\begin{array}{ll}{\lambda^{z},} & {z=0,1,2,3, \ldots} \\ {\mu^{|z|},} & {z=0,-1,-2,-3, \ldots}\end{array}\right.
$$
Then we have $\mathbb{E}|{{\tilde e}_i}| = \frac{{\lambda  + \mu }}{{1 - \lambda \mu }}$ and
$$
\operatorname{Var} Z=\frac{1}{(1-\lambda)^{2}(1-\mu)^{2}}\left(\frac{\mu(1-\lambda)^{3}(1+\mu)+\lambda(1-\mu)^{3}(1+\lambda)}{1-\lambda \mu}-(\lambda-\mu)^{2}\right)
$$
$$
\mathbb{P}(Z \leq x)=\left\{\begin{array}{ll}{\frac{(1-\lambda) \mu^{-[x]}}{1-\lambda \mu}} & {\text { if } x<0} \\ {1-\frac{(1-\mu) \lambda^{[x]+1}}{1-\lambda \mu}} & {\text { if } x \geq 0}\end{array}\right.
$$
Moreover, we have
\[P(|{{\tilde e}_i}| > c) = \frac{{(1 - \mu ){\lambda ^{[c] + 1}} + (1 - \lambda ){\mu ^{[c] + 1}}}}{{1 - \lambda \mu }}.\]
\[P(\mathop {\max }\limits_i |{{\tilde e}_i}| > c) = 1 - {\left( {1 - \frac{{(1 - \mu ){\lambda ^{[c] + 1}} + (1 - \lambda ){\mu ^{[c] + 1}}}}{{1 - \lambda \mu }}} \right)^n}.\]
{\bf Proof.} The formula of mean, absolute mean, variance, distribution function are a direct results from Proposition 2.2 of \cite{Kozubowski2006}.  From the expression of distribution functions, we get
\[P({{\tilde e}_i} > c) + P({{\tilde e}_i} <  - c) = \frac{{(1 - \mu ){\lambda ^{[c] + 1}}}}{{1 - \lambda \mu }} + \frac{{(1 - \lambda ){\mu ^{ - [ - c - 1]}}}}{{1 - \lambda \mu }} = \frac{{(1 - \mu ){\lambda ^{[c] + 1}} + (1 - \lambda ){\mu ^{[c] + 1}}}}{{1 - \lambda \mu }}.\]
From the property of maximum, we have
\[P(\mathop {\max }\limits_i |{{\tilde e}_i}| > c) = 1 - {[1 - P(|{{\tilde e}_i}| > c)]^n} = 1 - {\left( {1 - \frac{{(1 - \mu ){\lambda ^{[c] + 1}} + (1 - \lambda ){\mu ^{[c] + 1}}}}{{1 - \lambda \mu }}} \right)^n}.\]

 \end{lemma}

\section{The QQ plots with $\epsilon=\log(n)/n^{1/2}$}

\par This section contains two figures in the case of $\epsilon=\log(n)/n^{1/2}$. The first is the QQ-plots for $\hat{\xi}_{ij}$ with $q=3$, while the second is the QQ plots for $\hat{\xi}_{ij}$ and $\bar{\xi}_{ij}$ with $q=2$, where $q$ is the value of edge weights, $\hat{\xi}_{ij}$ and $\bar{\xi}_{ij}$ are the non-denoised and denoised estimates, respectively.

\begin{figure*}

\centering

\subfigure[$n=100,\epsilon=\log(n)/n^{1/2}$]{

\begin{minipage}[]{1\textwidth}
\centering
\includegraphics[width=0.44\textwidth]{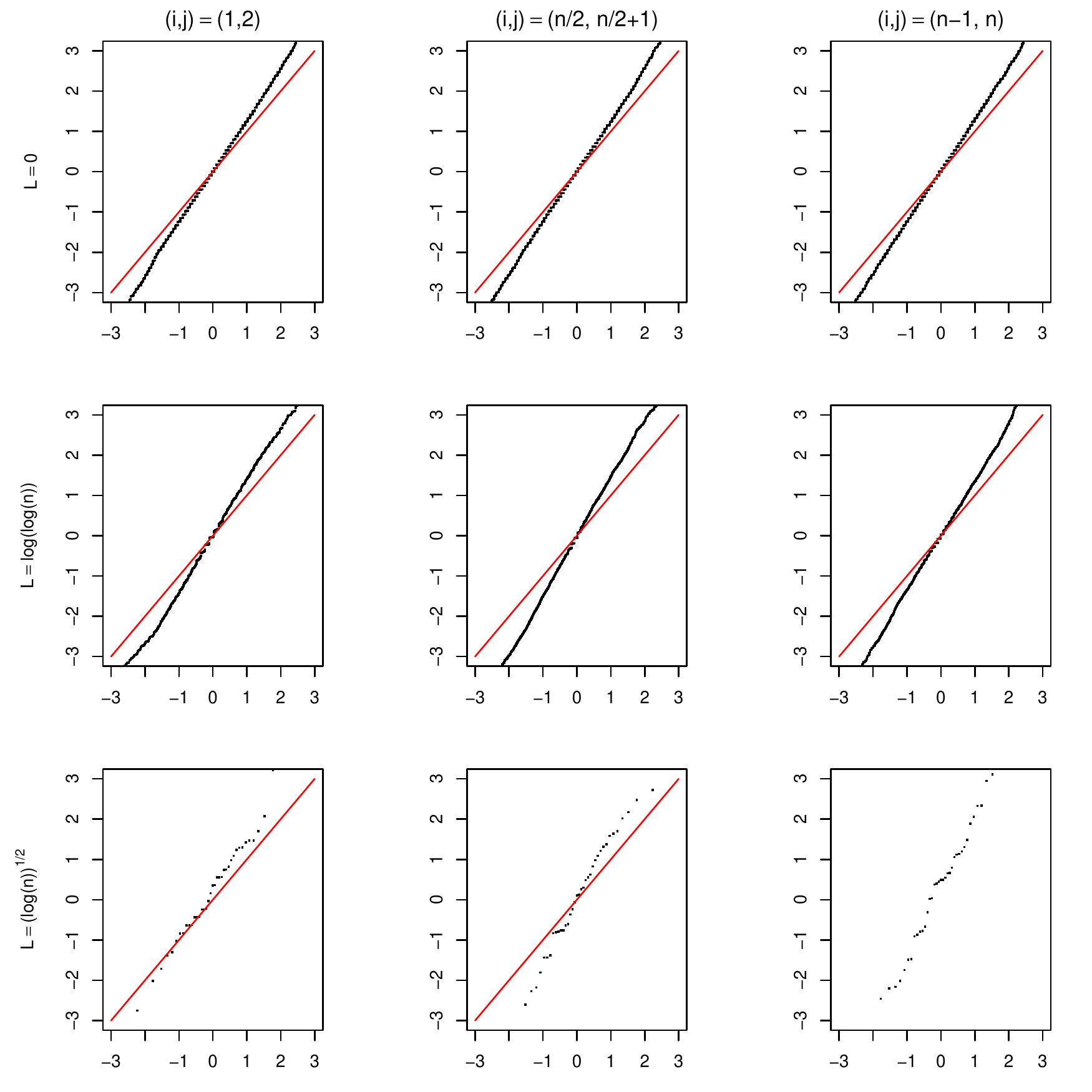}

\end{minipage}

}

\subfigure[$n=200,\epsilon=\log(n)/n^{1/2}$]{

\begin{minipage}[b]{1\textwidth}
\centering
\includegraphics[width=0.44\textwidth]{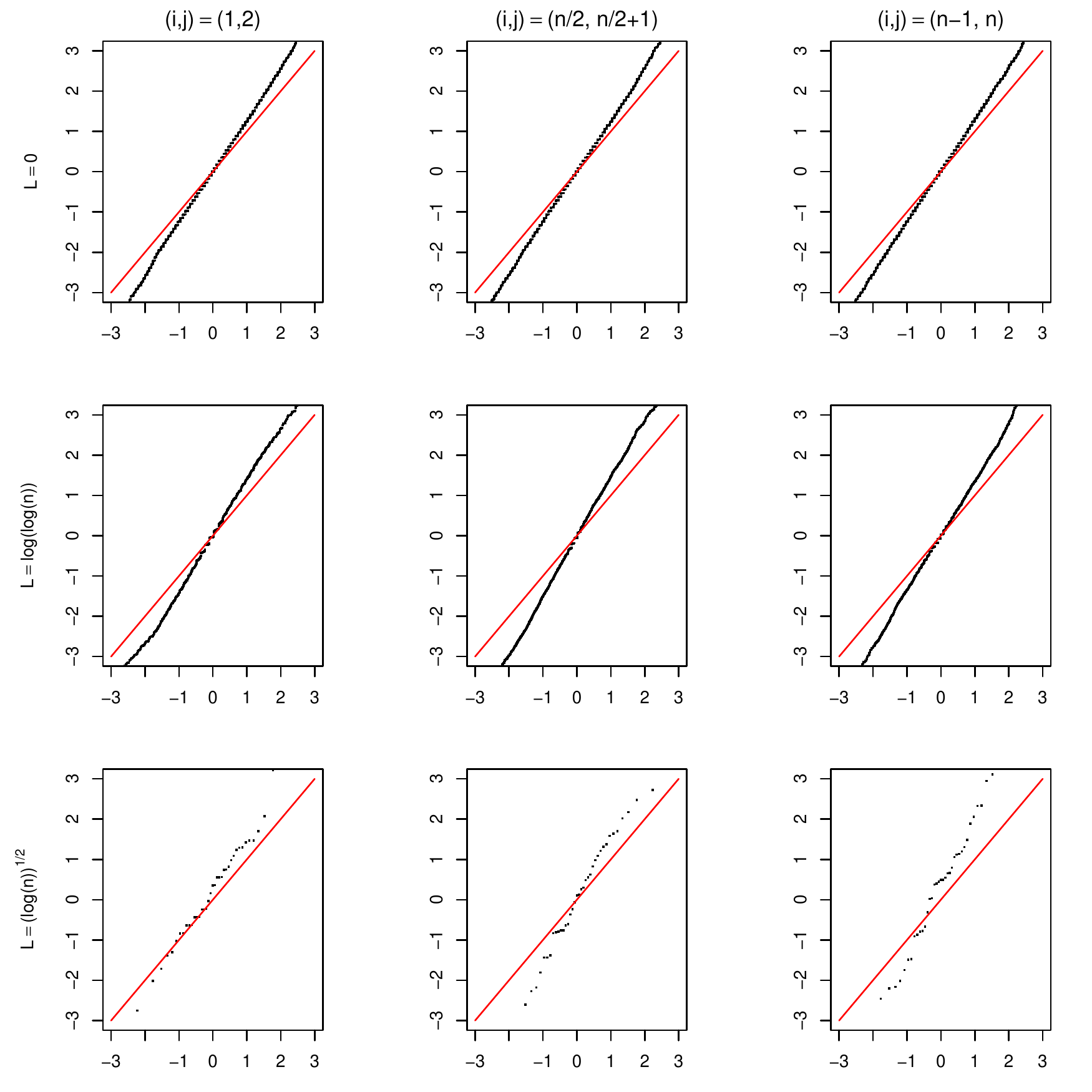}

\end{minipage}

}

\subfigure[$n=500,\epsilon=\log(n)/n^{1/2}$]{

\begin{minipage}[b]{1\textwidth}
\centering
\includegraphics[width=0.44\textwidth]{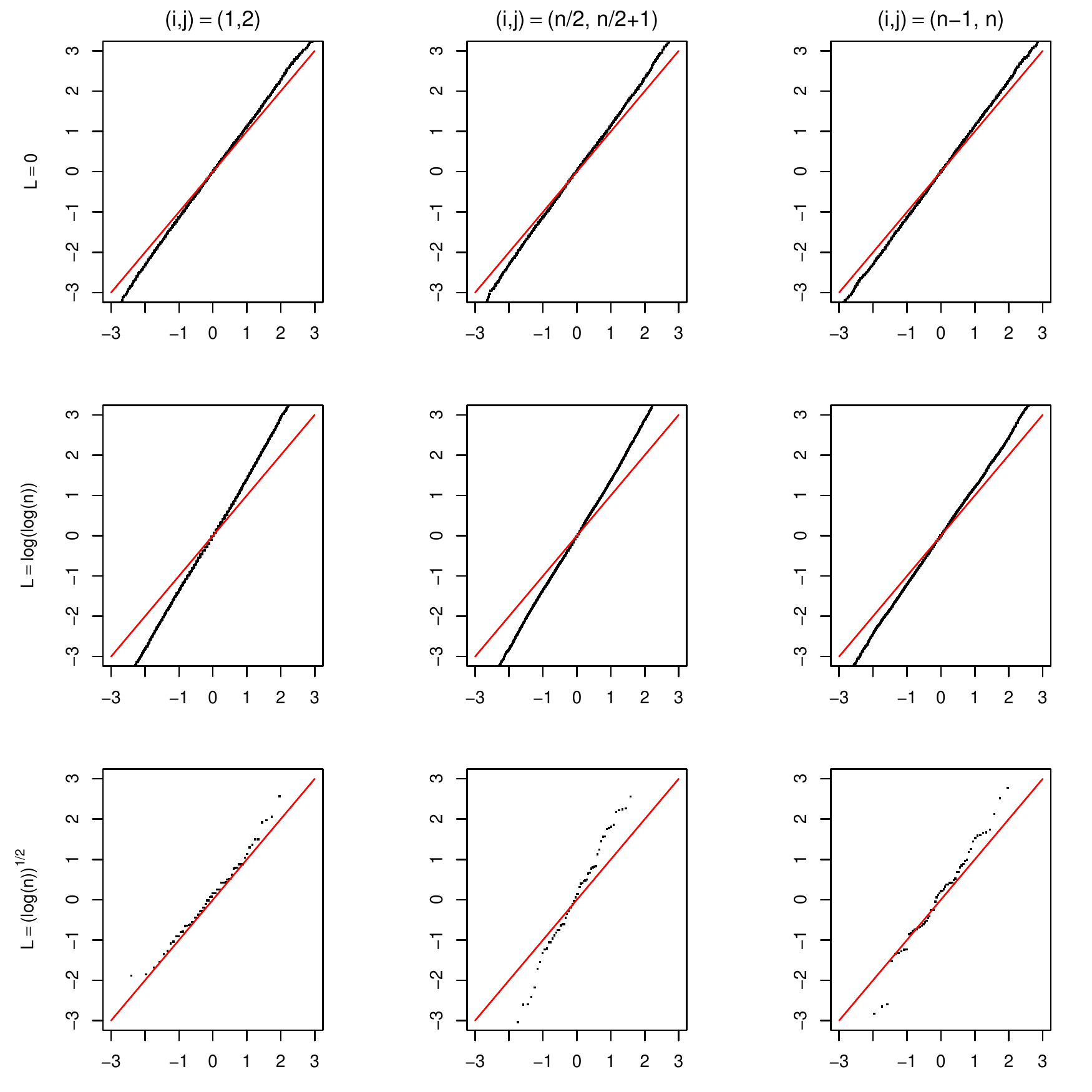}

\end{minipage}

}
\caption{The QQ plots with $q=3$}
\label{Afig1}
\end{figure*}

 \begin{figure}
\centering

\subfigure[$n=100,\epsilon=\log(n)/n^{1/2}$]{
\begin{minipage}[b]{1\textwidth}
\centering
\includegraphics[width=0.45\textwidth]{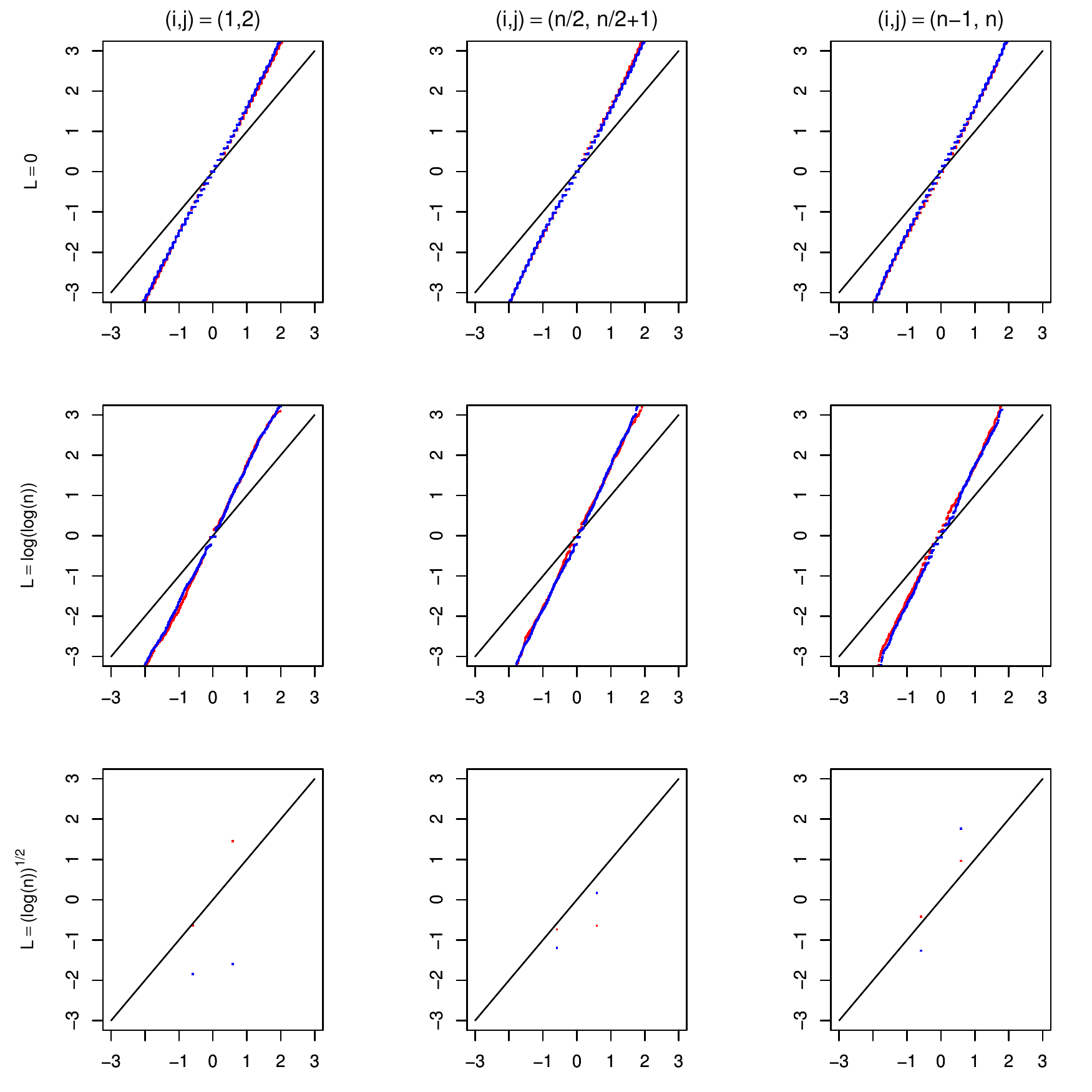}
\end{minipage}
}
\subfigure[$n=200,\epsilon=\log(n)/n^{1/2}$]{
\begin{minipage}[b]{1\textwidth}
\centering
\includegraphics[width=0.45\textwidth]{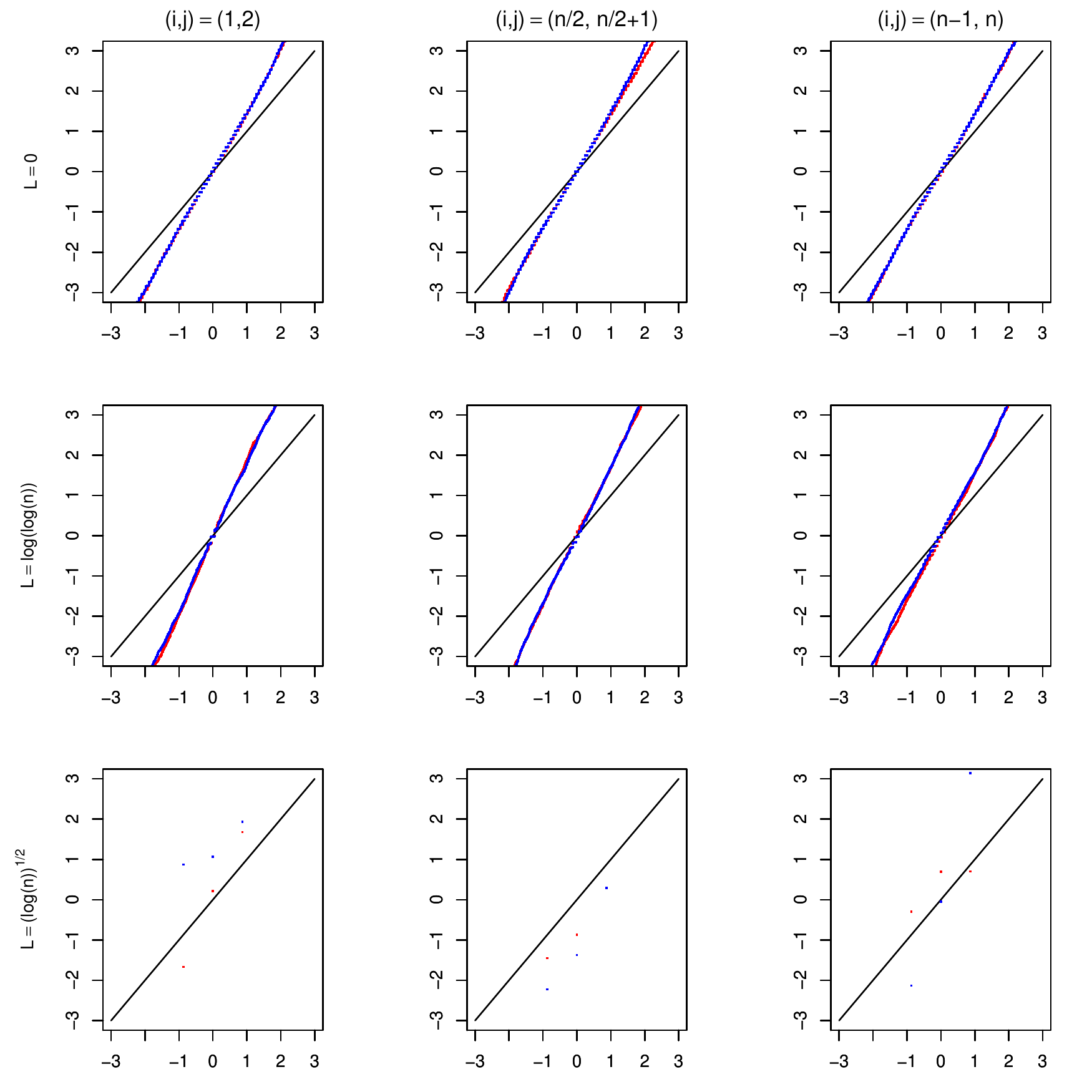}
\end{minipage}
}
\caption{The QQ plots of the non-denoised and denoised estimates with $q=2$}
\label{Afig2}
\end{figure}


\begin{thebibliography}{}
%%\bibliographystyle{imsart-nameyear}
\setlength{\itemsep}{-1.5pt}

\bibitem[Aggarwal and Yu(2007)]{Aggarwal and Yu 2007}
Aggarwal C. C., Yu P. S., 2007. On privacy-preservation of text and sparse binary data with sketches. Proceedings of the 2007 SIAM International Conference on Data Mining, 57--67.

\bibitem[Blitzstein and Diaconis(2011)]{Blitzstein and Diaconis(2011)}
Blitzstein J.,  Diaconis P., 2011. A sequential importance sampling algorithm for generating random graphs with prescribed degrees. Internet Math. 6, 489--522.

\bibitem[Chatterjee et al.(2011)]{Chatterjee:Diaconis:Sly:2011}
Chatterjee S., Diaconis P., Sly A., 2011. Random graphs with a given degree sequence. Ann. Appl. Probab. 21, 1400--1435.

\bibitem[Chung(2001)]{Chung2001}
Chung K. L., 2001. A course in probability theory, 3rd. Academic press.

\bibitem[Duncan et al.(2004)]{Duncan: Keller-McNulty:Stokes:2004}
Duncan G. T., Keller-McNulty S. A., Stokes S. L., 2004. Disclosure risk vs. data utility: the R-U confidentiality map. Technical Report Number 121
of National Institute of Statistical Sciences.

\bibitem[Duchi et al.(2018)]{Duchi2018}
Duchi J. C., Jordan M. I., Wainwright M. J., 2018. Minimax optimal procedures for locally private estimation. J. Am. Stat. Assoc.  113(521), 182--201.

\bibitem[Dwork(2006a)]{Dwork2006a}
Dwork C., 2006a. Differential privacy, in: Proceedings of the 33rd international colloquium on automata. Languages and Programming, 1--12.

\bibitem[Dwork et al.(2006b)]{Dwork:McSherry:Nissim:Smith:2006b}
Dwork C., McSherry F., Nissim K., Smith A, 2006b. Calibrating noise to sensitivity in private data analysis. In TCC, 265--284.


\bibitem[Farine(2014)]{Farine(2014)}
Farine D. R., 2014. Measuring phenotypic assortment in animal social networks: Weighted associations are more robust than binary edges. Anim. Behav. 89,
141--153.

\bibitem[Fung et al.(2007)]{Fung:Wang:Yu:2007}
Fung B. C. M., Wang K., Yu P. S, 2007. Anonymizing classification data for privacy preservation. IEEE Transactions on Knowledge and Data Engineering (TKDE) 19, 5(May), 711--725.

\bibitem[Ghinita et al.(2008)]{Ghinita:Tao:Kalnis:2008}
 Ghinita G., Tao Y., Kalnis P., 2008. On the anonymization of sparse high-dimensional data. In Proc.
of the 24th IEEE International Conference on Data Engineering (ICDE). 715--724.

\bibitem[Gragg and Tapia(1974)]{Gragg:Tapia:1974}
 Gragg W. B., Tapia  R. A., 1974. Optimal error bounds for the Newton-Kantorovich theorem. SIAM J. Numer. Anal. 11, 10--13.

\bibitem[Granovetter(1993)]{Granovetter(1993)}
 Granovetter M. S., 1993. The strength of weak ties. Am. J. Sociol. 1360--1380.

\bibitem[Haratym(2017)]{Granovetter(2017)}
Haratym E., 2017. Animals' right to privacy. World Scientific News, 85, 73--77.

\bibitem[Hay et al.(2009)]{Hay:Li:Miklau:Jensen:2009}
Hay M., Li C., Miklau G., Jensen D., 2009. Accurate estimation of the degree distribution of private networks. In Data Mining. ICDM'09. Ninth IEEE International Conference on 169--178.

\bibitem[Hillar et al.(2013)]{Hillar:Wibisono:2013}
Hillar C., Wibisono A., 2013.  Maximum entropy distributions on graphs. Available at \url{http://arxiv.org/abs/1301.3321.}

\bibitem[Holohan et al.(2017)]{Holohan:Leith:Mason:2017}
Holohan N., Leith D. J., Mason O., 2017. Extreme points of the local differential privacy polytope. Linear. Algebra. Appl. 534, 78--96.

\bibitem[Inusah and Kozubowski(2006)]{Inusah2006}
Inusah, S.,  Kozubowski, T. J. (2006). A discrete analogue of the Laplace distribution. J. Stat. Plan. Infrer. 136(3), 1090-1102.

\bibitem[Jackson.(2008)]{Jackson.(2008)}
Jackson M. O., 2008. Social and economic networks. Princeton University Press, Princeton.

\bibitem[Karwa et al.(2016)]{Karwa2016}
Karwa V., Slavkovi\'{c} A. B., 2016. Inference using noisy degrees: differentially private $\beta$-model and synthetic graphs. Ann. Stat. 44, 87--112.

\bibitem[Kasiviswanathan et al.(2013)]{Kasiviswanathan2013}
Kasiviswanathan S.P., Nissim K., Raskhodnikova S., Smith A. (2013). Analyzing Graphs with Node Differential Privacy. In: Sahai A. (eds) Theory of Cryptography. Lecture Notes in Computer Science, vol 7785. Springer, Berlin, Heidelberg.

\bibitem[Kotz et al.(2012)]{Kotz2012}
Kotz S., Kozubowski T., Podgorski K., 2012. The Laplace distribution and generalizations: a revisit with applications to communications, economics, engineering, and finance. Springer.

\bibitem[Kozubowski and Inusah(2006)]{Kozubowski2006}
Kozubowski T. J.,  Inusah S., 2006. A skew Laplace distribution on integers. Ann. I. Stat. Math. 58(3), 555--571.

\bibitem[Lang(1993)]{Lang1993}
Lang S., 1993. Real and Functional Analysis. Springer.

\bibitem[Lauritzen.(2008)]{Lauritzen.(2008)}
Lauritzen S. L., 2008. Exchangeable Rasch matrices. Rendiconti di Matematica, Series VII. 28, 83--95.

\bibitem[Li et al.(2007)]{Li:Li:Venkatasubramanian:2007}
Li N., Li T., Venkatasubramanian S., 2007.  $t-$closeness: Privacy beyond $k-$anonymity and $1-$diversity. Proceedings of the 23rd International Conference on Data Engineering 106--115.

\bibitem[Li et al.(2016)]{Li:Shen:Lang:Dong:2016}
Li Y., Shen H., Lang C., Dong H., 2016. Practical anonymity models on protecting private weighted graphs. Neurocomputing, 218, 359--370.

\bibitem[{Lounici(2008)}]{Lounici08}
Lounici, K. 2008. Sup-norm convergence rate and sign concentration property of Lasso and Dantzig estimators. Electron. J. Stat. 2, 90--102.

\bibitem[Machanavajjhala et al.(2006)]{Machanavajjhala: Gehrke:Kifer D:Venkitasubramaniam:2006}
 Machanavajjhala A., Gehrke J., Kifer D., Venkitasubramaniam M., 2006. $\mathscr{L}$-diversity: Privacy beyond kappa-anonymity. Proceedings of the 22nd International Conference on Data Engineering 24.

\bibitem[McSherry and Talwar(2007)]{McSherry :Talwar:2007}
 McSherry F., Talwar K., 2007. Mechanism design via differential privacy, in: Proceedings of the 48th Annual Symposium on Foundations of Computer Science. IEEE. 94--103.

\bibitem[Nissim et al.(2007)]{Nissim:Raskhodnikova:Smith:2007}
 Nissim K., Raskhodnikova S.,  Smith A., 2007. Smooth sensitivity and sampling in private data analysis. In Proceedings of the thirty-ninth annual ACM Symposium on Theory of Computing, ACM. 75--84.

\bibitem[Samarati and Sweeney(1998)]{Samarati and Sweeney1998}
Samarati P., Sweeney L., 1998. Protecting privacy when disclosing information: $k$-anonymity and its enforcement through generalization and suppression. Tech. rep., SRI International. March.

\bibitem[Steutel and van Harn(2003)]{Steutel2003}
Steutel F. W., van Harn K., 2003. Infinite divisibility of probability distributions on the real line. CRC Press.

\bibitem[Sundaresan et al.(2007)]{Sundaresan:Fischhoff:  Dushoff:  Rubenstein : 2007}
Sundaresan S. R., Fischhoff I. R., Dushoff J., Rubenstein D. I., 2007. Network metrics reveal diVerences in social organization between two Wssion-fusion species, Grevy's zebra and onager. Oecologia 151, 140--149.

\bibitem[{van der Vaart(1998)}]{Vaart1998}
Van der Vaart A. W., 1998. Asymptotic statistics (Vol. 3). Cambridge university press.

\bibitem[Wasserman and Zhou(2010)]{Wasserman:Zhou:2010}
Wasserman L., Zhou S., 2010. A Statistical Framework for Differential Privacy. J. Am. Stat. Assoc. 105(489), 375--389.

\bibitem[Yan and Xu(2013)]{Yan:Xu:2013}
Yan T., Xu J., 2013. A central limit theorem in the $\beta$-model for undirected random graphs with a diverging number of vertices. Biometrika 100, 519--524.

\bibitem[Yan et al.(2015)]{Yan:Zhao:Qin:2015}
Yan T., Zhao Y., Qin H., 2015. Asymptotic normality in the maximum entropy models on graphs with an increasing number of parameters. J. Multivariate Anal. 133, 61--76.

\bibitem[Yan and Leng(2015)]{Yan2015}
Yan T.,  Leng C., 2015. A simulation study of the $p_{1}$ model for directed random graphs. Stat. Interface, 8, 255--266.

\bibitem[Yan et al.(2016a)]{Yan:Leng:Zhu:2016a}
Yan T., Leng C., Zhu J., 2016a. Asymptotics in directed exponential random graph models with an increasing bi-degree sequence. Ann. Stat. 44, 31--57.

\bibitem[Yan et al.(2016b)]{Yan:Qin:Wang:2016b} Yan T., Qin H., Wang H., 2016b. Asymptotics in undirected random graph models parameterized by the strenghs of vertices. Stat. Sinica. 26, 273--293.

\bibitem[{Zhang et al.(2014)}]{Zhang14}
Zhang H., Liu Y., Li B. 2014. Notes on discrete compound Poisson model with applications to risk theory. Insur. Math. Econ. 59, 325--336.

\end{thebibliography}

\newpage
\begin{thebibliography}{}
%%\bibliographystyle{imsart-nameyear}
\setlength{\itemsep}{-1.5pt}
\bibitem[Fan et al.(2020)]{AFan2020}
Fan Y., Zhang H., Yan T., 2020. Asymptotic Theory for Differentially Private Generalized $\beta$-models with Parameters Increasing. Submitted.
\bibitem[Karwa et al.(2016)]{AKarwa2016}
Karwa V., Slavkovi\'{c} A. B., 2016. Inference using noisy degrees: differentially private $\beta$-model and synthetic graphs. Ann. Stat. 44, 87-112.

\bibitem[Kotz et al.(2012)]{AKotz2012}
Kotz S., Kozubowski T., Podgorski K., 2012. The Laplace distribution and generalizations: a revisit with applications to communications, economics, engineering, and finance. Springer.

\bibitem[Kozubowski and Inusah(2006)]{AKozubowski2006}
Kozubowski T. J.,  Inusah S., 2006. A skew Laplace distribution on integers. Ann. I. Stat. Math. 58(3), 555-571.

\end{thebibliography}
\end{document}